

\magnification=\magstep1
\documentstyle{amsppt}
\NoBlackBoxes
\def\bh{\Cal B(\Cal H)}
\def\on{\Cal O_n}
\def\h{\Cal H}
\def\1{\hbox{$1\hskip -3.75pt 1$}}
\def\End{\operatorname{End}}
\def\exp{\operatorname{exp}}
\def\aut{\operatorname{Aut}}
\def\shift{\operatorname{Shift}}
\def\uhf{\operatorname{UHF}}
\def\irr{\operatorname{Irr}}
\def\rep{\operatorname{Rep}}
\def\erg{\operatorname{Erg}}
\def\ad{\operatorname{Ad}}
\def\ind{\operatorname{Index}}
\def\Dim{\operatorname{Dim}}
\def\linspan{\operatorname{lin\,span}}
\def\tr{\operatorname{tr}}
\def\id{\operatorname{id}}
\topmatter
\title Endomorphisms of $\bh$ \endtitle
\author Ola Bratteli, Palle E.T. Jorgensen, and Geoffrey L.
Price \endauthor
\address Mathematics Institute, University of Oslo, PB 1053
-- Blindern, N-0316 Oslo, Norway \endaddress
\email bratteli\@math.uio.no \endemail
\address Department of Mathematics, University of Iowa, Iowa
City, IA  52242, U.S.A. \endaddress
\email jorgen\@math.uiowa.edu\endemail
\address Department of Mathematics 9E, U.S. Naval Academy,
Annapolis, MD  21402, U.S.A. \endaddress
\email GLP\@sma.usna.navy.mil \endemail
\thanks Research supported in part by the Norwegian Research
Council for Ola Bratteli and Palle Jorgensen, by the U.S.
National Science Foundation for Palle Jorgensen, and by the
U.S. National Security Agency for Geoffrey Price. \endthanks
\keywords Endomorphisms, shift, Cuntz-algebra, pure states,
infinite tensor product, Hilbert space \endkeywords
\subjclass 46L10, 46L50, 47A58, 47C15, 81S99 \endsubjclass
\abstract The unital endomorphisms of $\bh$ of (Powers)
index $n$ are classified by certain $U(n)$-orbits in the set
of non-degenerate representations of the Cuntz algebra $\on$
on $\h$. Using this, the corresponding conjugacy classes are
identified, and a set of labels is given. This set of labels
is $P/\sim$ where $P$ is a set of pure states on the $\uhf$-
algebra $M_{n^\infty}$, and $\sim$ is a non-smooth
equivalence on $P$. Several subsets of $P$, giving concrete
examples of non-conjugate shifts, are worked out in detail,
including sets of product states, and a set of nearest
neighbor states. \endabstract
\endtopmatter

\document

\head 0. Introduction\endhead
Recently the study of endomorphisms of von Neumann algebras
has received increased attention, both in connection with
the Jones index for subfactors and its applications
\cite{Jon}, and also in connection with duality for compact
groups \cite{Wor} and super-selection sectors in algebraic
quantum field theory. Two other articles (by W.~Arveson and
by R.~Powers) in these proceedings deal with {\it
semigroups\/} of endomorphisms of the type I$_\infty$-
factor. Here we restrict to the case of {\it single\/}
endomorphisms of $\bh$. Potentially it is expected that the
theory for $\bh$ may possibly be extended or modified to
apply also to other factors, but so far only a few
relatively isolated results (although still some very
important ones) are known for endomorphisms of factors other
than $\bh$. We report here on recent and new developments in
the study of $\End(\bh)$. The methods used draw among other
things on seminal ideas of von Neumann, and also on ideas of
Powers from his pioneering work on the states on the CAR
(canonical anticommutation relation)-algebra, and, more
generally, states on the UHF (uniformly hyperfinite) $C^*$-
algebras.

The work on $\End(M)$ for the case when $M$ is a von Neumann
factor of type II$_1$ (especially the hyperfinite case) is
ongoing. It will not be treated here, but we refer to
\cite{Pow2}, \cite{Po-Pr}, \cite{EW}, \cite{Cho}, and
\cite{ENWY}.

\head 1. Main Results\endhead
Let $\bh$ be the $C^*$-algebra of bounded linear operators
on the separable, infinite dimensional Hilbert space $\h$.
If $\alpha:\bh\rightarrow\bh$ is a unital endomorphism, we
say that $\alpha$ is ergodic if $\{X\in \bh \mid
\alpha(X)=X\}=\Bbb C 1$, and that $\alpha$ is a shift if
$\bigcap_{n=1}^\infty \alpha^n (\bh)=\Bbb C 1$. The (Powers)
index $n\in \{1,2,\ldots,\infty\}$ of $\alpha$ is defined as
the $n$ such that $\alpha(\bh)' \cap \bh$ is isomorphic to
the factor of type I$_n$, \cite{Pow2}. Thus $n=1$ if and
only if $\alpha$ is an automorphism. Let $\End_n(\bh)$
(respectively $\erg_n(\bh)$, $\shift_n(\bh)$) denote the set
of unital endomorphisms (respectively ergodic unital
endomorphisms, shifts) of $\bh$ of index $n$. We say that
two elements $\alpha,\beta\in\End(\bh)$ are conjugate if
there is an automorphism $\gamma\in\aut(\bh)=\End_1(\bh)$
such that $\alpha=\gamma \circ \beta\circ \gamma^{-1}$, and
$\alpha,\beta$ are approximately conjugate if for all
$\epsilon>0$ there is a $\gamma\in \aut(\bh)$ such that
$\|\alpha - \gamma\circ\beta\circ\gamma^{-1} \|<\epsilon$.
It is easy to see that any two approximately conjugate
endomorphisms $\alpha,\beta$ must have the same index $n$.

In \cite{Pow2, Theorem 2.3} it was proved that if
$\alpha,\beta$ are shifts of index $n\geq 2$ each allowing a
pure, normal invariant state on $\bh$, then $\alpha$ and
$\beta$ are conjugate. The problem of whether there exist
shifts without invariant vector states was left open in
\cite{Pow2}, but we will both construct explicit classes of
examples of shifts of order $n$ without invariant vector-
states in Sections 5--8, and prove a classification theorem.

Our construction of these special shift-conjugacy classes,
and our analysis of their ergodic theoretic, and clustering
type properties, are based on fundamental ideas of von
Neumann, especially his 1938 {\it Compositio-paper\/}
\cite{vNeu}, and their extension by Guichardet \cite{Gui}
(notably \cite{Gui3}). The imprint on our paper from von
Neumann's legacy is most visible in our construction of
explicit examples in Sections 6, 7, and 8 below.

In the study of $\End(\bh)$ we will make extensive use of
ideas developed by von Neumann and other pioneers in
operator algebras and in quantum theory, \cite{vNeu},
\cite{Seg1--2}, \cite{Pow1}, \cite{ArWoo} (see also the
beginning of Remarks 8.2).

\proclaim{Theorem 1.1} \rom{(see \cite{Lac1, Theorem 4.5})}
Assume $n\in\{2,3,4,\ldots,\infty\}$. Then the set of
conjugacy classes in $\shift_n (\bh)$ can be equipped with a
natural Borel structure which is not countably separated.
The same applies to $\End_n(\bh)$ and $\erg_n(\bh)$. In
particular there exist elements in $\shift_n(\bh)$ which do
not allow invariant vector states.\endproclaim

This theorem will be proved in Section 5 (the Borel
structure portion is new). In Section 5 we will give a
complete labeling of the conjugacy classes in
$\shift_n(\bh)$ by $P/ \sim$, where $P$ is a subset of the
pure state space of the $\uhf$ algebra $M_{n^\infty}$, and
$\sim$ is a certain equivalence relation on $P$. In Sections
5, 6, 7, and 8  we will look at some special elements in
$P/\sim$.  On the way to the proof, we will gain further
insight into the shifts allowing invariant vector states.

In \cite{Lac1}, M. Laca continues the program initiated by
Powers of analyzing the conjugacy classes of discrete shifts
on $\bh$. The central theme of his approach, as it is here,
is to exploit the correspondence between endomorphisms and
representations of the Cuntz algebras which implement the
endomorphisms. In his paper Laca succeeds in establishing
the existence of uncountably many conjugacy classes of
shifts of each index \cite{Lac1, Remark 4.6.2}. He also
obtains \cite{Lac1, Theorem 4.5} a characterization of the
conjugacy classes of shifts which identifies them with an
equivalence class structure of a certain family of pure
states on the subalgebra $\uhf_n$ of the Cuntz algebra
$\on$. This result appears in a slightly different guise as
our Theorem 1.1, which is included for the purposes of
exposition. In \cite{Pow2, Theorem 2.4} it was proved that
any two shifts of $\bh$ with the same index are outer
conjugate. Another version of this result is:

\proclaim{Theorem 1.2} \rom{(see \cite{Pow2} and \cite{Lac1,
Proposition 2.3})} Let $\alpha,\beta$ be two endomorphisms
of $\bh$ of the same index $n\in\{1,2,\ldots,\infty\}$. Then
there is a unitary $U\in\bh$ such that
$$
\alpha(X)=U\beta(X)U^*
$$
for all $X\in\bh$.\endproclaim

Defining $\gamma (X)=UXU^*$, this relation can also be
expressed as
$$\align
\alpha(X)&=U\beta(U^*UXU^*U)U^*\\
   &=\gamma\beta\gamma^{-1}(UXU^*)
\endalign$$
which is the form of outer conjugacy considered in
\cite{Pow2}. We will see that one cannot in general find a
unitary $U$ such that $\alpha(X)=\beta (UXU^*).$ This is
proved in Section 3. Finally, using Voiculescu's non-
commutative Weyl--von Neumann theorem \cite{Voi1, Wor}, we
can establish

\proclaim{Theorem 1.3} Let $\alpha,\beta$ be two
endomorphisms of $\bh$ of the same index $n$, $2\leq n
<\infty$. Then $\alpha$ and $\beta$ are approximately
conjugate; i.e., there is a sequence $\gamma_k \in
\aut(\bh)$ such that
$$
\| \alpha - \gamma_k \circ \beta \circ \gamma_k^{-
1}\|\rightarrow 0.
$$
The sequence $\gamma_k$ may furthermore be chosen such that
$\alpha(X)- (\gamma_k\circ\beta\circ\gamma_k^{-1})(X)$ is
compact for each $X\in\bh$, $k\in \Bbb N$.\endproclaim

We remark that when $n=1$, it is well known that an
automorphism $\alpha$ of $\bh$ is implemented by a unitary
operator $U$, unique up to a scalar phase factor, and thus
$\aut(\bh)$ is indexed by the set $\rep(C(\Bbb T),\h)$ of
non-degenerate representations of $C(\Bbb T)$ on $\h$,
modulo the canonical action of the circle group $\Bbb T$.
These representations are well known from spectral theory,
\cite{Ped}. Thus $\shift_1(\bh)$ and $\erg_1(\bh)$ are
empty, and $\End_1(\bh)$ is countably separated in its
natural Borel structure. Theorem 1.2 is trivially true in
the case $n=1$ (just put $U=U_\alpha U_\beta^*$ where
$\alpha=\ad(U_\alpha)$, $\beta=\ad(U_\beta))$, while Theorem
1.3 is false.

This work was essentially completed before we became aware
of Laca's results. As mentioned above, some of our work
overlaps with that in \cite{Lac1}, and we indicate below
where this occurs. Our approach to the subject differs in
several aspects, however. A major goal of our work, for
example, is to develop techniques and concepts which
differentiate between those endomorphisms which admit normal
invariant states and those which do not (all endomorphisms
have invariant states, however, see Remark 7.6). Since
Powers already showed that for each index there is only one
conjugacy class of {\it shifts\/} allowing invariant normal
pure states (see Theorem 4.2, below), any method giving
other conjugacy classes of course gives shifts without
invariant vector states. (There does, however, exist a
plethora of conjugacy classes of non-ergodic endomorphisms
of a given index $n$ each allowing (several) invariant
vector states; just take discrete direct sums of the
representations of $\on$ defined by Cuntz's states as in
Section 4.) A special feature of our approach is that we
obtain many examples of shifts not allowing invariant vector
states by perturbing shifts allowing such states by various
perturbation techniques (see Sections 6 and 7). Our
constructions in Sections 5--7 are based primarily on
consideration of (infinite) {\it product states\/} on
$\uhf_n$, whereas our construction in Section 8 uses instead
certain {\it nearest neighbor states\/} on $\uhf_n$. Both
constructions lead to shifts which do not have invariant
vector states, but, more importantly, the shifts on nearest
neighbor states are {\it not\/} conjugate to those from
Sections 6--7.

In Section 9, we construct explicitly extensions of
endomorphisms of $\bh$ to automorphisms of $\Cal B(\h
\otimes\h)$.

We finally point out the connection between our results and
the results of Arveson on one-parameter semigroups of $^*$-
endomorphisms (see \cite{Arv1--2}). If one translates
Arveson's concepts, which are tailor-made for the semigroup
$\Bbb R_+$, to the semigroup $\Bbb N \cup \{0\}$, then his
{\it spectral $C^*$-algebra\/} for a shift of index $n$ is
nothing but the Toeplitz-Cuntz algebra $\Cal E_n$, which in
turn is an extension of $\on$ by the compact operators when
$n$ is finite, and $\Cal E_\infty = \Cal O_\infty$
\cite{Eva}. Otherwise Arveson's Fock space methods have a
different flavor from our infinite tensor product methods.

The Toeplitz-Cuntz algebras also play a role in the recent
work in Dinh \cite{Din}, as well as \cite{Lac1--2} and
\cite{Sta}.

\head 2. Cuntz Algebras and Cuntz States\endhead
The Cuntz algebra $\on$ is uniquely defined as the $C^*$-
algebra generated by $n=2,3,\ldots$ isometries
$s_1,\ldots,s_n$ satisfying
$$
s_i^*s_j=\delta_{ij}1,\qquad \sum_{j=1}^n s_j s_j^* =1, \tag
2.1
$$
\cite{Cun}. There is a canonical representation of the $n$-
dimensional unitary group $U(n)$ in the automorphism group
of $\on$ defined by
$$
\tau_g(s_i)=\sum_{j=1}^n g_{ji}s_j\tag 2.2
$$
for $g=[g_{ij}]_{i,j=1}^n \in U(n)$.

Let $\pi_1,\pi_2$ be two non-degenerate representations of
$\on$ on a Hilbert space $\h$, and put
$$
S_i=\pi_1(s_i),\qquad T_i=\pi_2(s_i).\tag 2.3
$$
Then there exists a unitary operator $M=[m_{ij}]\in M_n
(\bh)$ and a  unitary operator $U\in\bh$ such that
$$
T_i=\sum_{j=1}^n S_j m_{ji}=US_i.\tag 2.4
$$
The operators $M$ and $U$ are given uniquely by
$$
m_{ji}=S_j^* T_i,\qquad U=\sum_{j=1}^n T_j S_j^*\tag 2.5
$$
and we have the relations
$$
m_{ji}=S_j^* US_i,\qquad U=\sum_{j,i=1}^n S_j m_{ji}
S_i^*.\tag 2.6
$$
Conversely, if $\{S_i\}$ is a realization of $\on$ on $\h$,
and $[m_{ij}]$ is a unitary element in $M_n(\bh)$, then
$\{T_i\}$ defined by \rom{(2.4)} is a realization of $\on$
on $\h$. The same remark applies to a single unitary
operator $U\in \bh$, and the other equation in \rom{(2.4)}.
We will give explicit formulas for the transfer operators
\rom{(2.6)} in Sections 7--8 below for elements in
$\End_n(\bh)$ from distinct conjugacy classes.

The $C^*$-algebra $\on$ is simple and antiliminal when
$n>1$, \cite{Cun}. We define, naturally, $\Cal O_1$ as the
universal $C^*$-algebra generated by one unitary element,
i.e., $\Cal O_1 =C(\Bbb T)$, and $\Cal O_\infty$ as the
algebra generated by isometries $s_1,s_2,\ldots$ satisfying
merely the first relation in \rom{(2.1)}. Then $\Cal
O_\infty$ is still simple and antiliminal \cite{Cun}, while
$\Cal O_1$ of course is abelian.

With a slight abuse of terminology, we will say that $\pi$
is a non-degenerate representation of $\Cal O_\infty$ if
$\pi$ is a representation with $\sum_{i=1}^\infty \pi (s_i
s_i^*)=1$, where the sum is in the strong operator topology.
With this convention, all the statements in the second
paragraph of this section are still valid for $n=\infty$,
and the infinite sums converge in the strong operator
topology.

Let $\uhf_n$ be the fixed point subalgebra of $\on$ under
the canonical action of the center of $U(n)$. Thus $\uhf_n$
is the closure of the linear span of operators of the form
$$
s_{i_1}s_{i_2}\cdots s_{i_k}s_{j_k}^* s_{j_{k-1}}^* \cdots
s_{j_1}^*
$$
over $k=0,1,2,\ldots$. If $n<\infty$, then $\uhf_n$ is the
$\uhf$-algebra $M_{n^\infty}$, which is the uniform closure
of finite linear combinations of operators of the form $A_1
\otimes A_2 \otimes A_3 \otimes \cdots$, where each $A_i$
acts on a fixed $n$-dimensional Hilbert space (i.e., $A_i\in
M_n$) and all but finitely many of the $A_i$'s are the
identity. If $n=\infty$, then $\uhf_\infty$ is the (non-
simple) AF algebra described as follows: Let $\h$ be a fixed
infinite-dimensional separable Hilbert space. For each $k\in
\Bbb N$, let $\Cal C_k$ be the $C^*$-algebra of compact
operators on $\bigotimes_1^k \h$, viewed as the $C^*$-
algebra generated by all linear combinations of elements of
the form $A_1 \otimes A_2 \otimes \cdots \otimes A_k \otimes
I \otimes I\otimes \cdots$, where $A_i \in C(\h)$. Then
$\uhf_\infty$ is the $C^*$-algebra generated by the $\Cal
C_k$'s for $k\in \Bbb N$, and the identity. For more details
on $\uhf_\infty$, see also \cite{Cun}, \cite{Eva}, \cite{Br-
Rob, Example 5.3.27} or \cite{Lac1--2}.

Let $D_n$ denote the canonical diagonal subalgebra of
$\uhf_n$; that is, $D_n$ is the abelian $C^*$-algebra
obtained as the closure of the linear span of
$$
s_{i_1} s_{i_2} \cdots s_{i_k} s_{i_k}^* s_{i_{k-1}}^*
\cdots s_{i_1}^*.
$$
Then $D_n$ is maximal abelian in $\uhf_n$. If $2\leq n <
\infty$ then $D_n$ is canonically isomorphic to $C(\prod
_{k=0}^\infty \Bbb Z_n)$, where $\Bbb Z_n =\{1,\ldots,n\}$
equipped with the discrete topology. If $n=\infty$, $D_n$ is
the abelian $C^*$-algebra spanned by $1\otimes 1\otimes
1\otimes \cdots, C_0(\Bbb N)\otimes 1\otimes 1\otimes\cdots,
C_0(\Bbb N \times \Bbb N)\otimes 1\otimes\cdots$. (For
details on this, see \cite{Br-Rob; Example 5.3.27}.)

If $n<\infty$, we defined the canonical endomorphism $\psi$
of $\on$ by
$$
\psi(x)=\sum_{k=1}^n s_k x s_k^*.
$$
Then $\left. \psi\right|_{\uhf_n}$ is the one-sided shift.

If $\eta_1, \ldots ,\eta_n\in \Bbb C$ with $\sum_{i=1}^n
|\eta_i |^2 =1$ the associated Cuntz state is the pure state
$\omega_\eta$ on $\on$ defined by
$$
\omega_\eta (s_{i_1}\cdots s_{i_k}s_{j_1}^* \cdots
s_{j_\ell}^*)
=\eta_{i_1}\cdots \eta_{i_k}\bar\eta_{j_1}\cdots
\bar\eta_{j_\ell}
$$
(this definition also goes through with obvious
modifications for $n=\infty$ and $n=1$). When $2\leq n
<\infty$, $\left. \omega_\eta \right|_{\uhf_n} $ is the
infinite product on $\bigotimes_0^\infty M_n$ of the pure
states on $M_n$ defined by the vector $\eta=(\eta_1,
\ldots,\eta_n)$. When $n=+\infty$, $\left.
\omega_\eta\right|_{\uhf_\infty}$ is similarly the state on
$\uhf_\infty$, described as before, defined by the unit
vector $\eta\otimes\eta\otimes\eta \otimes \cdots$,
\cite{Cun}, \cite{ACE}, \cite{BEGJ} and \cite{Voi2}.

\head 3. Endomorphisms\endhead
\proclaim{Theorem 3.1} \rom{(\cite{Arv1}, \cite{Lac1;
Theorem 2.1, Proposition 2.2})} Let $\varphi$ be a unital
endomorphism of $\bh$ of Powers index $n\in
\{1,2,3,\ldots,+\infty\}$.

It follows that there exists a non-degenerate representation
of $\on$ on $\h$ such that
$$
\varphi(X)=\sum_{i=1}^n S_i X S_i^*\tag 3.1
$$
where $S_i$ is the representative of $s_i$. Conversely, any
non-degenerate representation of $\on$ on $\h$ defines an
endomorphism of index $n$ by \rom{(3.1)}. The representation
is unique up to the canonical action of $U(n)$.\endproclaim

\demo{Proof} Since $\varphi(\bh)$ is a unital subalgebra of
$\bh$, isomorphic to $\bh$, we have a tensor product
decomposition $\h=\h_0 \otimes \Cal K$ of $\h$ such that
$\varphi(\bh)$ identifies with $\Cal B(\Cal H_0)\otimes 1$,
and then $\varphi (\bh)'\cap \bh \cong 1\otimes \Cal B(\Cal
K)$, \cite{Dix}. Thus, $\ind(\varphi)=\Dim(\Cal K)$.

Let $(E_{ij})_{i,j=1}^n$ be a complete set of matrix units
for $\varphi(\bh)'$. It follows that
$$
E_{ii}\varphi(\bh)=\varphi(\bh)E_{ii}\cong \Cal B(\h_0)\cong
\bh
$$
for $i=1,2,\ldots,n$, and $X\rightarrow \varphi(X)E_{ii}$ is
a $*$-isomorphism between $\bh$ and $\Cal B(E_{ii}\h)$. By
Wigner's theorem (which is Theorem 3.1 in the case $n=1$)
there is a unitary operator $S_i$ from $\h$ onto $E_{ii}\h$
such that
$$
\varphi(X)E_{ii}=S_i X S_i^*.
$$
But then
$$
\varphi(X) =\varphi(X)\sum_{i=1}^n E_{ii}
   =\sum_{i=1}^n \varphi(X)E_{ii}
  =\sum_{i=1}^n S_i X S_i^*.
$$
We have
$$
S_i^* S_i =1,\qquad \sum_{i=1}^n S_iS_i^*=\sum_{i=1}^n
E_{ii}=1
$$
so the $S_i$ satisfy the Cuntz relations \rom{(2.1)}.
Conversely, if $S_i$ satisfy the Cuntz relations, then
$\varphi$ defined by \rom{(3.1)} is an endomorphism such
that $\varphi(\bh)' \cap \bh$ is spanned by $S_iS_j^*$, and
consequently $\varphi(\bh)'\cap\bh\cong M_n$ and $\varphi$
has index $n$.

Let $T_i$, $i=1,\ldots,n$ be another non-degenerate
realization of $\on$ that implements~$\varphi$:
$$
\varphi(X)=\sum_{i=1}^n T_i X T_i^* = \sum_{i=1}^n S_i X
S_i^*.
$$
Multiply the last relation to the left by $S_j^*$ and to the
right by $T_i$ to obtain
$$
S_j^* T_i X = X S_j^* T_i.
$$
Since this is true for any $X\in \bh$,
$$
S_j^* T_i =h_{ji}1
$$
where $h_{ji}\in \Bbb C$. But then $h=[h_{ji}]\in U(n)$ and
$$
\pi_2 =\pi_1 \circ \tau_h
$$
where $\tau$ is the canonical action of $U(n)$ on $\on$, and
$\pi_1$, $\pi_2$ are the representations determined by $S$,
$T$, respectively.\enddemo

\definition{Definition 3.2} For $n=1,2,\ldots,\infty$, let
$$
\rep(\on,\h)
$$
denote the set of all non-degenerate representations of
$\on$ on $\h$, and
$$
\irr (\on,\h)
$$
the set of all irreducible representations of $\on$ on $\h$,
and
$$
\rep_s(\on,\h)
$$
the set of representations of $\on$ on $\h$ such that
$\uhf_n$ is weakly dense in $\bh$. Of course the two latter
sets are empty if $n=1$.\enddefinition

The canonical action of $U(n)$ on $\on$ gives rise to an
action of $U(n)$ on each of these spaces. Also, the unitary
group $U(\h)$ on $\h$ acts on each of the three spaces by
$\pi(\cdot) \rightarrow U\pi (\cdot)U^*$ for $U\in U(\h)$,
$\pi\in\rep(\on,\h)$. The following  corollary of Theorem
3.1 is immediate.

\proclaim{Theorem 3.3} Let $\pi \rightarrow\varphi(\pi)$ be
the surjective map from $\rep(\on,\h)$ onto \linebreak
$\End_n(\bh)$ defined in Theorem 3.1. Then:\roster

\item"{(3.2)}" $\varphi(\pi) \in\erg_n (\bh)$ if and only if
$\pi\in\irr (\on,\h)$

\item"{(3.3)}" $\varphi(\pi) \in \shift_n(\bh)$ if and only
if $\pi\in\rep_s (\on,\h)$

\item"{(3.4)}" \rom{(\cite{Lac, Proposition 2.4})}
$\varphi(\pi_1)$ and $ \varphi(\pi_2)$ are conjugate if and
only if there is a $g\in U(n)$ and a $U\in U(\h)$ such that
$$
\pi_2(\cdot)=U\pi_1 (\tau_g(\cdot))U^*.
$$
\endroster

\noindent In short, the conjugacy classes in $\End_n(\bh)$
correspond to the orbits in \linebreak $\rep(\on,\h)$ under
the joint actions of $U(n)$ and $U(\h)$.\endproclaim

\demo{Proof} To prove \rom{(3.2)}, it suffices to show that
(\cite{Lac, Proposition 3.1})
$$
\pi(\on)'=\{X\in \bh\mid\varphi(\pi)(X)=X\}\equiv
\bh^\varphi. \tag 3.5
$$
But if $X\in\pi(\on)'$, then
$$
\varphi(X)=\sum_{i=1}^n S_iX S_i^*=\sum_{i=1}^n
S_iS_i^*X=1\cdot X=X
$$
where $S_i=\pi(s_i)$, so $X\in\bh^\varphi$, and
$\pi(\on)'\subseteq \bh^\varphi$. Conversely, if
$X\in\bh^\varphi$, then $\sum_{i=1}^n S_i X S_i^* =X$.
Multiplying to the left by $S_j^*$ we obtain
$$
XS_j^* =S_j^*X
$$
and since
$X^*\in\bh^\varphi$, we also derive
$$
S_jX=XS_j.
$$
Hence $X\in\pi(\on)'$ and so
$$
\bh^\varphi \subseteq \pi(\on)'.
$$
This establishes \rom{(3.5)} and therefore \rom{(3.2)}.

To prove \rom{(3.3)} we will more generally establish that
(\cite{Lac, Proposition 3.1})
$$
\bigcap_k \varphi^k (\bh)= (\pi(\uhf_n))'\cap \bh.\tag 3.6
$$
This again will follow from
$$
\varphi^k (\bh)' =\linspan \{S_{i_1}\cdots S_{i_k} S_{j_k}^*
\cdots S_{j_1}\}.
\tag 3.7
$$
But as
$$
\varphi^k(X)=\sum_{i_1,\ldots, i_k=1}^n
S_{i_1} \cdots S_{i_k} X S_{i_k}^* \cdots S_{i}^*
$$
and $(i_1 , \ldots, i_k) \rightarrow S_{i_1} \cdots S_{i_k}$
is a non-degenerate representation of $\Cal O_{n^k}$, it
suffices to prove \rom{(3.7)} for $k=1$. But as
$$\align
S_iS_j^* \varphi(X) &= S_iS_j^* \sum_k S_k X S_k^*\\
   &= S_j X S_j^* = \sum_k S_k X S_k^* S_i S_j^*\\
   &= \varphi(X) S_i S_j^*
\endalign$$
we have
$$
\linspan \{S_i S_j^*\} \subseteq \varphi(\bh)'.
$$
Conversely, a general element $X\in\bh$ may be written
$$
X=\sum_{ij} S_i X_{ij} S_j^*
\qquad\text{where}\qquad
X_{ij}=S_i^* X S_j
$$
and, if $X\in \varphi (\bh)'$, then
$$
\sum_{ij} S_i X_{ij} S_j^*  \sum_k S_k Y S_k^*
=\sum_k S_k Y S_k^* \sum_{ij} S_i X_{ij} S_j^*
$$
for all $Y\in\bh$; that is,
$$
\sum_{ij} S_i X_{ij} Y S_j^* = \sum_{ij} S_i YX_{ij} S_j^*
$$
for all $Y\in\bh$. Thus $X_{ij}$ must be scalar multiples of
$1$, and $X$ is a linear combination of $S_i S_j^*$. This
establishes \rom{(3.7)}, and hence \rom{(3.6)} and
\rom{(3.3)}. (Of course, if $n=\infty$, linear span means
the weak closure of the linear span.)

To prove \rom{(3.4)}, put $S_i=\pi_1 (s_i)$, $T_i=\pi_2
(s_i)$. If $\varphi(\pi_1)$, $\varphi(\pi_2)$ are conjugate,
there exists a $\gamma=\ad(U) \in \aut(\bh)$ such that
$$
\varphi(\pi_2)=\gamma \varphi (\pi_1)\gamma^{-1}
$$
i.e.,
$$
\sum_i T_i X T_i^* = U\left(\sum_{i} S_i U^* X U S_i^*
\right) U^*
    = \sum_i (US_i U^*) X (US_i U^* )^*
$$
for all $X\in\bh$. Since $US_i U^*$ satisfy the Cuntz
relations it follows from the uniqueness part of Theorem 3.1
that there exists a $g=[g_{ij}]\in U(n)$ such that
$$
T_i = \sum_{j=1}^n g_{ji} US_j U^*,
$$
i.e.,
$$
\pi_2(\cdot) =U(\pi_1 \circ \tau_g(\cdot))U^*.
$$

The converse is established by doing the steps in converse
order.

This finishes the proof of Theorem 3.3.\qed\enddemo

\demo{Proof of Theorem 1.2} By Theorem 3.1 there exist two
realizations $S,T$ of $\on$ on $\h$ such that
$$
\alpha(X) =\sum_i S_i X S_i^*,\qquad \beta(X)=\sum_i T_i X
T_i^*.
$$
By \rom{(2.4)} there is a unitary $U$ such that
$$
S_i =UT_i
$$
for all $i$. But then
$$
\alpha(X)=U\beta (X)U^*
$$
for all $X\in \bh$.\qed
\enddemo

\demo{Proof of Theorem 1.3} By Theorem 3.1 there exist two
realizations $S,T$ of $\on$ on $\h$ such that
$$
\alpha(X) = \sum_{i=1}^n S_i X S_i^*,\qquad
\beta(X)=\sum_{i=1}^n T_i X T_i^*.
$$
Let $\epsilon >0$. As $\on$ is a simple, antiliminal $C^*$-
algebra it follows from Voiculescu's non-commutative Weyl--
von Neumann theorem (\cite{Voi1, Corollary 1.4} and
\cite{Wor}) that there exists a unitary $U$ on $\h$ such
that
$$
S_i - U T_i U^*
$$
are compact for $i=1,\ldots,n$, and
$$
\|S_i -UT_i U^* \| <\epsilon/2n.
$$
Let $\gamma(X) =U X U^*$. Then
$$\align
\alpha(X) - \gamma \beta \gamma^{-1} (X)
  &= \sum_{i=1}^n \left( S_i X S_i^* -UT_iU^* X
(UT_iU^*)^*\right)\\
  &= \sum_{i=1}^n (S_i - UT_i U^*) X S_i^*
       +\sum_{i=1}^n UT_i U^* X (S_i-UT_i U^*)^*.
\endalign$$
Thus $\alpha(X)-\gamma\beta\gamma^{-1}(X)$ is compact and
$$\align
\|\alpha(X) - \gamma\beta\gamma^{-1} (X)\|
  & \leq 2n \cdot 1 \cdot \|X\| \epsilon/2n\\
  & = \epsilon\|X\|.
\endalign$$
This proves Theorem 1.3.\qed\enddemo

\head 4. Shifts and Invariant States \endhead
Let $\alpha$ be an endomorphism of $\bh$. The next theorem
gives a characterization of the normal $\alpha$-invariant
pure states on $\bh$.

\proclaim{Theorem 4.1} Let $\alpha$ be a unital endomorphism
of $\bh$ of index $n = 1,2, \ldots,\infty$, and let $\pi$ be
a corresponding non-degenerate representation of $\on$. Let
$S_i =\pi(s_i)$, $i=1,\ldots,n$. Let $\xi$ be a unit vector
in $\h$, and let $\omega(X) = \langle \xi, \pi(X)\xi
\rangle$ be the corresponding state on $\on$. The following
conditions are equivalent:
\roster
\item"{(4.1)}" $\langle \xi,\alpha(X)\xi\rangle=\langle
\xi,X\xi \rangle$ for all $X\in\bh$.
\item"{(4.2)}" $\xi$ is a joint eigenvector for $S_i^*$ for
$i=1,2,\ldots,n$.
\item"{(4.3)}" $\omega$ is a Cuntz state on $\on$.
\endroster
Furthermore, the corresponding eignvalues in \rom{(4.2)} are
$\bar\eta_i$:
$$
S_i^*\xi=\bar\eta_i\xi \tag 4.4
$$
for $i=1,2,\ldots,n$, if and only if $\sum_{i=1}^n
|\eta_i|^2=1$ and $\omega = \omega_\eta$.\endproclaim

\demo{Proof} (4.2) $\Leftrightarrow$ (4.3) and the final
remark are straightforward.

(4.2) $\Rightarrow$ (4.1): If $S_i^* \xi = \bar\eta_i \xi$,
then
$$\align
\sum_{i=1}^n |\eta_i|^2 &=\sum_{i=1}^n \langle \bar\eta_i
\xi, \bar\eta_i \xi \rangle
  =\sum_{i=1}^n \langle S_i^* \xi, S_i^* \xi \rangle\\
 &=\sum_{i=1}^n \langle \xi, S_iS_i^* \xi\rangle
  =\langle \xi, \xi \rangle =1.
\endalign$$
Furthermore
$$\langle \xi, \alpha(X)\xi \rangle
  =\sum_{i=1}^n\langle S_i^* \xi, XS_i^* \xi \rangle
  = \sum _{i=1}^n |\eta_i|^2 \langle \xi,X\xi \rangle
  =\langle \xi, X\xi \rangle.
$$

(4.1) $\Rightarrow$ (4.2): Assume that $\langle \xi,\alpha
(X)\xi\rangle = \langle \xi, X\xi \rangle$. We have
$$
\langle \xi, \alpha(X)\xi \rangle = \sum_{i=1}^n \langle
S_i^* \xi, X S_i^* \xi \rangle.
$$
But $X \rightarrow \langle S_i^*\xi, X S_i^*\xi\rangle$ is a
positive linear functional on $\bh$ of norm
$$
\langle S_i^*\xi,S_i^* \xi \rangle= \langle \xi, S_iS_i^*\xi
\rangle.
$$
The sum of these norms is
$$
\sum_{i=1}^{n} \langle \xi, S_i S_i^* \xi \rangle
  =\left\langle \xi, \left( \sum_{i=1}^n S_iS_i^*\right)\xi
\right\rangle
  = \langle \xi,\xi  \rangle =1.
$$
As the sum of these positive functionals is $\langle \xi,
\cdot \xi\rangle$ and $\langle \xi, \cdot \xi\rangle$ is
pure, it follows that
$$
\langle S_i^*\xi, X S_i^* \xi\rangle
=\|S_i^* \xi \|^2 \langle \xi, X\xi\rangle
$$
for all $X\in \bh$, but then
$$
S_i^* \xi = \bar\eta_i\xi
$$
where $\eta_i \in \Bbb C$ is such that $|\eta_i| = \|S_i^*
\xi\|$.\qed\enddemo

Using Theorem 4.1 we can prove the following result, which
is implicit in the proof of Theorem 2.3 of \cite{Pow2}, see
also \cite{Sta} for related results.

\proclaim{Theorem 4.2} Suppose that $\alpha$, $\beta$ are
ergodic unital endomorphisms of $\bh$, both of index $n\in
\{2,3,\ldots \}$ and assume that both $\alpha$ and $\beta$
allow a pure invariant state.

It follows that $\alpha$ and $\beta$ are conjugate, and both
of them are shifts.\endproclaim

\demo{Proof} Let $\pi_\alpha$, $\pi_\beta$ be the
representations of $\on$ corresponding to $\alpha$, $\beta$,
respectively. The ergodicity of $\alpha$, $\beta$ implies
that $\pi_\alpha$, $\pi_\beta$ are irreducible, by Theorem
3.3. Let $\xi_\alpha$, $\xi_\beta$ be unit vectors in $\h$
such that $\langle \xi_\alpha,\alpha (\cdot)\xi_\alpha
\rangle = \langle \xi_\alpha, \cdot \xi_\alpha \rangle$ and
$\langle \xi_\beta, \beta(\cdot) \xi_\beta \rangle = \langle
\xi_\beta, \cdot \xi_\beta \rangle$. By Theorem 4.1 the
corresponding two states on $\on$ are Cuntz states; i.e.,
there exist unit vectors $\eta^\alpha =
(\eta_1^\alpha,\ldots, \eta_n^\alpha)$ and $\eta^\beta =
(\eta_1^\beta, \ldots, \eta_n^\beta)$ in
$\ell^2(\{1,2,\ldots, n\})$ such that
$$
\langle \xi_\alpha, \pi_\alpha (x)\xi_\alpha \rangle
   =\omega_{\eta^\alpha}(x)\qquad\text{and}\qquad
\langle \xi_\beta, \pi_\beta (x)\xi_\beta\rangle
  =\omega_{\eta^\beta} (x)
$$
for $x\in\on$. Now, choose $g=[g_{ij}]\in U(n)$ so that
$\eta^\alpha = g^T \eta^\beta$, where $g^T$ is the transpose
of $g$. But since
$$
\omega_\eta (\tau_g (s_i))
  =\omega_\eta \left(\sum_j g_{ji}s_j\right)
  =\sum_j g_{ji} \eta_j
  =\omega_{g^T \eta}(s_i)
$$
etc., one has
$$
\omega_\eta \circ \tau_g = \omega_{g^T \eta}
$$
for any $g\in U(n)$ and any unit vector $\eta\in \ell^2
(\{1,2,\ldots\})$. In particular
$$
\omega_{\eta^\beta} \circ \tau_g = \omega_{\eta^\alpha}.
\tag 4.5
$$
As $\pi_\alpha$ and $\pi_\beta$ are irreducible, it follows
that $\xi_\alpha$, is cyclic for $\pi_\alpha$, and
$\xi_\beta$ is cyclic for $\pi_\beta$, and hence the
relation \rom{(4.5)} entails that $\pi_\alpha$ and
$\pi_\beta \circ \tau_g$ are unitarily equivalent. By
\rom{(3.4)}, $\alpha$ and $\beta$ are conjugate. To show
that $\alpha$ and $\beta$ are shifts is equivalent to
showing that $\pi_\alpha(\uhf_n)$ and $\pi_\beta (\uhf_n)$
are weakly dense in $\bh$.  But $\pi_\alpha$, $\pi_\beta$
are unitarily equivalent to the representation defined by
the Cuntz's states $\omega_{\eta^\alpha}$,
$\omega_{\eta^\beta}$, and these are irreducible in
restriction to $\uhf_n$ by \rom{(8.8)--(8.9)}
below.\qed\enddemo

\head 5. Classification of Conjugacy Classes of Shifts
\endhead

In this section we will prove Theorem 1.1, and find an
explicit set of labels for the conjugacy classes in
$\shift_n(\bh)$. In Sections 6, 7, and 8 we will consider
some more explicit points in the label space. Assume that
$n\in \{2,3,\ldots\}$.  The case $n=\infty$ is somewhat more
complicated and was treated in detail in \cite{Lac}. The
results are similar in that case, and we will restrict to
finite $n$ in the rest of this section. Consider unital
shifts of Powers index $n$ on $\bh$. By Theorem 3.3, these
correspond to the set $\rep_s(\on,\h)$ of representations
$\pi$ of $\on$ on $\h$ such that $\pi(\uhf_n)$ is weakly
dense in $\bh$. These representations identify with the
cyclic representation defined by any vector state, defined
by a unit vector in $\h$. We will characterize abstractly
the corresponding states on $\on$, or, rather, the
restriction of those states to $\uhf_n$. So let $P$ denote
the set of pure states $\omega$ on $\uhf_n$ such that
$\omega$ has a pure extension $\omega'$ to $\on$ with the
property that, if
$(\h_{\omega'},\pi_{\omega'},\Omega_{\omega'})$ is the
corresponding representation, then
$\pi_{\omega'}(\uhf_n)''=\Cal B(\h_{\omega'})$. Let:
$$\align
\sigma(\cdot)&=\sum_iS_i\cdot S_i^* \text{ be the canonical
endomorphism of $\uhf_n$}\\
&\text{($=$ the one-sided shift on $M_{n^\infty}$)}\\
A_m &= \undersetbrace m \to{M_n \otimes \cdots \otimes M_n}
\otimes 1 \otimes 1 \otimes \cdots \subseteq \uhf_n \\
A_m^c &= \undersetbrace m \to{1\otimes \cdots \otimes 1}
\otimes M_n \otimes M_n \otimes \cdots \subseteq
  \uhf_n\\
&=\text{relative commutant of $A_m$}.
\endalign$$

Then $\sigma(A^c_m)\subseteq A_{m+1}^c$ and $\sigma:A_m^c
\rightarrow A_{m+1}^c$ is an isomorphism.

\proclaim{Lemma 5.1}  If $\omega$ is a pure state on
$\uhf_n$ then $\omega \circ \sigma$ is a type I factor state
of multiplicity $\leq n$.\endproclaim

\demo{Proof}
$$
\pi_\omega (\sigma(\uhf_n))' = \pi_\omega(A_1^c)\cong M_n,
$$
and the representation Hilbert space of $\omega\circ\sigma$
identifies with $\overline{\pi_\omega
(\sigma(\uhf_n))\Omega_\omega}$.
\enddemo

\proclaim{Lemma 5.2} Let $\omega$ be a pure state on
$\uhf_n$. The following conditions are equivalent:
\roster\item"{(1)}" $\omega\in P$
\item"{(2)}" For all $\epsilon >0$ there is an $m\in \Bbb N$
such that
$$
\left\| \left. (\omega \circ \sigma -
\omega)\right|_{A_m^c}\right\| < \epsilon
$$
\item"{(3)}" $\lim_{m \rightarrow \infty} \|\omega \circ
\sigma^{m+1} - \omega\circ \sigma^m\| =0$\endroster
\endproclaim

\demo{Proof}Since $\sigma^m$ maps $\uhf_n$ isometrically
onto $A_m^c$, the equivalence of (2) and (3) is immediate.
Since $\omega$ and $\omega\circ\sigma$ both are factor
states by Lemma 5.1, it follows from \cite{Pow1, Theorem
2.7} that (2) is equivalent to the representations
$\pi_\omega$ and $\pi_{\omega\circ\sigma}$ being quasi-
equivalent.

(1) $\Rightarrow$ (2). If $\omega\in P$, then
$$
\omega \circ\sigma (x) =\sum_{i=1}^n \langle S_i^*
\Omega_\omega, \pi_\omega(x)S_i^* \Omega_\omega\rangle
$$
for $x\in \uhf_n$, where $S_i^*$ are the representatives of
$s_i^*$ in the extension of $\pi_\omega$ to a representation
of $\on$ on $\h_\omega$. But this shows that $\omega\circ
\sigma$ is a normal state in the representation
$\pi_\omega$, and, as $\omega$ and $\omega \circ \sigma$ are
factor states, they are quasi-equivalent.

(2) $\Rightarrow$ (1). If $\omega$ and $\omega\circ\sigma$
are quasi-equivalent, then the endomorphism $\pi_\omega (x)
\rightarrow \pi_\omega(\sigma(x))$, $x\in\uhf_n$, extends by
continuity to an endomorphism of $\Cal B(\h_\omega)$ which
we also call $\sigma$. But as $\pi_\omega (A_1) \subseteq
\pi_\omega (\sigma(\uhf_n))'$, we have $\pi_\omega (A_1)
\subseteq \sigma(\Cal B (\h_\omega))'$. Realizing the
elements in $\uhf_n$ as $n\times n$ matrices with entries in
$A_1^c$, using that $A_1 \cong M_n$, one easily deduces the
converse implication, and hence $\sigma$ has Powers index
$n$, and there exists by Theorem 3.1 a non-degenerate
representation $\pi$ of $\on$ on $\h_\omega$ such that
$$
\sigma(X) = \sum_{i=1}^n S_i X S_i^*
$$
where $S_i =\pi(s_i)$. But then $\sigma(\bh)'$ is spanned
linearly by $S_iS_j^*$, $i, j= 1\cdots n$, and, as
$\sigma(\bh)'=\pi_\omega(A_1)$, $\pi_\omega(A_1)$ is just
the linear span of $S_i S_j^*$, $i,j=1\cdots n$. Now,
modifying the $S_i$'s with an element in $U(n)$ if
necessary, we may arrange it so that
$$
S_iS_j^* =\pi_\omega (s_is_j^*)
$$
and this determines the $S_i$'s up to a fixed phase factor.
If
$$
e_{ij}^{(k)} = \sigma^k (s_i s_j^*)\qquad\text{and}\qquad
E_{ij}^{(k)} = \sigma^k (S_i S_j^*)
$$
then
$$
E_{ij}^{(k)} = \sigma ^k (\pi _\omega (s_i s_j^*))
=\pi_\omega (\sigma^k (s_is_j^*))
       =\pi_\omega(e_{ij}^{(k)})
$$
for $k=1,2,\ldots$, and thus we see that $\pi_\omega$
extends to a representation $\pi$ of $\on$ by setting
$$\pi(s_i)=S_i.
$$
Thus $\omega\in P$. \qed\enddemo

\proclaim{Lemma 5.3} Two elements $\omega,\omega' \in P$
define unitarily equivalent representations of\/ $\uhf_n$,
if and only, for $\forall \epsilon >0$, $\exists\, m$ such
that
$$
\|\left.( \omega - \omega')\right|_{A_m^c} \|<\epsilon.
$$
\endproclaim

\demo{Proof} \cite{Pow1, Theorem 2.7} again.
\enddemo

\proclaim{Lemma 5.4} Assume that $\omega,\omega' \in P$. The
following conditions are equivalent:
\roster\item"{(1)}" $\omega$ and $\omega'$ define conjugate
endomorphisms of $\bh$.
\item"{(2)}" There is a $g\in U(n)$ such that, for all
$\epsilon >0$, there is an $m\in \Bbb N$ with
$$
\left \|\left. (\omega-\omega' \circ \tau_g)\right|_{A_m^c}
\right \| < \epsilon
$$
where $\tau_g= \bigotimes_{k=1}^\infty \ad g$.
\item"{(3)}" There is a $g\in U(n)$ such that
$$
\lim_{m\rightarrow \infty} \| \omega \circ \sigma^m -
\omega' \circ \tau_g \circ \sigma^m \| =0.
$$
\item"{(4)}" There is a $g\in U(n)$ and a unitary $U\in
\uhf_n$ such that
$$
\omega(\cdot) = \omega' (U \tau_g(\cdot) U^*).
$$\endroster\endproclaim

\demo{Proof} The equivalence of the three first conditions
follows from Lemma 5.3 and Theorem 3.3. Since condition (2)
means that the two pure states $\omega$ and $\omega' \circ
\tau_g$ define unitarily equivalent representations,
condition (4) follows from Kadison's transitivity theorem,
\cite{KR}, and conversely (4) implies that $\omega$ and
$\omega' \circ \tau_g$ define unitary equivalent
representations. The proof is completed.\qed\enddemo

We are now ready to prove Theorem 1.1, and to even give an
explicit parametrization of the conjugacy classes in
$\shift_n(\bh)$. Let as before $P$ be the set of pure states
on $\uhf_n$ such that
$$
\lim_{m\rightarrow \infty} \| \omega \circ \sigma^{m+1}
- \omega \circ \sigma^m\| =0
$$
(this characterization is equivalent to the one given
above). Define two states $\omega,\omega' \in P$ to be
equivalent, $\omega \sim \omega'$, if they lie in the same
orbit in $P$ under the joint action of $U(n)$, and of
$U(\uhf_n)=$ the unitary group of $\uhf_n$. Then it follows
from Lemma 5.2, Lemma 5.4, and Theorem 3.3, that there is a
bijection between $P/\sim$ and the set of conjugacy classes
of endomorphisms of $\bh$. Since $\on$ is type III, and
$U(n)$ is compact, it follows, by the same reasoning as in
Glimm's theorem (see \cite{Gli} and \cite{Ped}), that
$P/\sim$ is not a standard Borel space. This is also implied
by the fact that the orbits in $\End_n(\bh)$ under conjugacy
all are norm dense by Theorem 1.3.

The slightly different proof in the case $n=\infty$ can be
found in \cite{Lac2}.

\example{Example 5.5} \rom{(\cite{Lac1})} Let $\xi_m,\xi_m'$
be unit vectors in $\Bbb C^n$, and let $\omega_m = \langle
\xi_m,\cdot \xi_m \rangle$ and $\omega_m'= \langle
\xi_m',\cdot \xi_m' \rangle$ be the corresponding pure
states on $M_n$. Consider the infinite tensor product states
$\omega = \bigotimes_{m=1}^\infty \omega_m$ and $\omega' =
\bigotimes_{m=1}^\infty \omega_{m}'$ on $\uhf_n =
\bigotimes_{m=1}^\infty M_n$. These are pure states, and by
Lemma 5.3 they induce unitarily equivalent representations
if and only if
$$
\lim_{m\rightarrow \infty} \|(\omega - \omega')\circ
\sigma^m \| =0 .\tag 5.1
$$
It is well known, \cite{Gui}, that this condition can be
expressed in the following equivalent ways:
$$\align
\sum_{m=1}^\infty \| \omega_n - \omega_n ' \|^2  &<\infty,
\tag 5.2\\
\sum_{m=1}^\infty (1 -  |\langle \xi_m, \xi_m' \rangle | )
&<\infty \tag 5.3\\
\lim_{k \rightarrow \infty}  \prod_{m=k}^\infty | \langle
\xi_m,\xi_m' \rangle | &=1.\tag 5.4
\endalign$$
These conditions are non-commutative versions of the
conditions for equivalence of infinite product measures on
$\prod_1^\infty \Bbb Z_n$ given by Kakutani in 1948,
\cite{Kak}. Similar conditions for quasi-equivalence of
quasi-free states, which are closely related to product
states, have been given in \cite{Po-St}, \cite{Ara1},
\cite{Ara2}, \cite{Dae}.

If, furthermore, the phases of the vectors $\xi_m'$ are
chosen optimally with respect to $\xi_m$, i.e., such that
$\langle \xi_m,\xi_m' \rangle \in \Bbb R_+$, then
\rom{(5.2)} is equivalent to
$$
\sum_{m=1}^\infty \| \xi_m - \xi_m' \|^2 < \infty. \tag 5.5
$$
Note for example that the equivalence of \rom{(5.5)} and
\rom{(5.3)} follows from
$$
\| \xi_m - \xi_m' \|^2 = 2(1-\operatorname{Re} \langle
\xi_m,\xi_m' \rangle).
$$
Using this, and Lemma 5.2, we see that $\omega =
\bigotimes_{m=1}^\infty \omega_m$ is in $P$ if and only if
$$
\sum_{m=1}^\infty \|\omega_m - \omega_{m+1} \|^2 <\infty;
\tag 5.6
$$
or, equivalently
$$
\sum_{m=1}^\infty \left(1- |\langle \xi_m, \xi_{m+1} \rangle
| \right) <\infty,
$$
or,
$$
\lim_{k \rightarrow \infty} \prod_{m=k}^\infty | \langle
\xi_m, \xi_{m+1} \rangle | =1,
$$
or, if the phases of $\xi_m$ are chosen inductively such
that $\langle \xi_m, \xi_{m+1} \rangle\in \Bbb R_+$,
$$
\sum_{m=1}^\infty \| \xi_m - \xi_{m+1} \| ^2 <\infty.
$$
The \rom{(5.6)} conditions are taken up again in Lemma 6.5
below. In Section 6 we will consider a condition \rom{(6.2)}
which is stronger than \rom{(5.6)}.

Finally, assume that $\omega = \bigotimes_{m=1}^\infty
\omega_m$ and $\omega' =\bigotimes_{m=1}^\infty \omega_m'$
are both in $P$. By Lemma 5.4, and the remarks above,
$\omega$ and $\omega'$ define non-conjugate shifts if and
only if, for all $g\in U(n)$
$$
\sum_{m=1}^\infty \| \omega_m - \omega_m' \circ \tau_g \| ^2
= +\infty \tag 5.7
$$
or, equivalently
$$
\sum_{m=1}^\infty \left(1- |\langle \xi_m, g\xi_m' \rangle |
\right) =+ \infty;
$$
or the other two similar conditions. In this way we may
analyze equivalence classes among the product states in $P$.
See Section 6 for more details.\endexample

\example{Example 5.6} Another way of constructing a
continuum of nonconjugate shifts is the following: Let
$(\lambda_i,\omega_i)_{i=1}^k $ be a finite sequence where
$\omega_i$ are distinct pure states on $M_n$, $\lambda_i >0$
and $\sum_{i=1}^k \lambda_i =1$. Choose $m$ so large that
$m^n \geq k$, and define a state $\omega$ on $A_m^c$ by
$$
\omega = \sum_{i=1}^k \lambda_i\, \undersetbrace {m+1 \text{
to }\infty}
\to{\omega_i \otimes\omega_i \otimes \cdots} . \tag 5.8
$$
By a standard construction (see \cite{Bra} and \cite{Gli}),
$\omega$ has an extension to a pure state on $\uhf_n$. This
extension $\omega$ has the property that
$\omega\circ\sigma^{j+1} = \omega \circ \sigma^j$ for $j
\geq m$, and hence $\omega \in P$. But it follows from Lemma
5.4 that a given pair, $(\lambda_i,\omega_i)_{i=1}^k$ and
$(\lambda_i',\omega_i')_{i=1}^{k'}$, gives rise to conjugate
shifts if and only if $k=k'$, and there exists a permutation
$\varphi$ of $\{1,\ldots,k\}$ and a $g\in U(n)$ such that
$$
\lambda_i =\lambda_{\varphi(i)}'
$$
and
$$
\omega_i = \omega_{\varphi(i)}' \circ \tau_g
$$
for $i=1,\ldots,k$. If $k=1$, this gives rise only to the
one conjugacy class allowing invariant vector states, but
for $k=2,3,\ldots$ there is a continuum  of distinct
possibilities.\endexample

\example{Example 5.7} As stated in Theorem 1.1, and
clarified in Theorem 1.3 and the remarks after Lemma 5.4,
there does not exist a smooth labeling of all the conjugacy
classes in $\shift_n(\bh)$, although there are of course
subclasses with a smooth labeling like those described in
Example 5.6. We will now give a complete labeling of a class
of shifts which will be described in more detail in Section
7, but again this labeling cannot be taken to be smooth.
Actually the conjugacy classes of the shifts obtained in
this fashion agree exactly with those obtained in Example
5.5, and these classes contain the classes described in more
detail in Sections 6 and 7 as subclasses. Let $e_i$,
$i=0,\ldots,n-1$ be the orthonormal basis of $\Bbb C^n$
defined by \rom{(7.14)}, and define $\bigotimes_{m=0}^\infty
\Bbb C^n \simeq \Cal L^2 (\Omega, \mu)$ as in the
introduction to Section 7, so that $\mu$ is normalized Haar
measure on $\Omega =\prod_{m=0}^\infty \Bbb Z_n$. In
particular $\Cal L^2(\Omega,\mu)$ contains the vector
$$
\1 = \bigotimes_{m=0}^\infty e_0 = e_0 \otimes e_0 \otimes
\cdots
$$\endexample

The following result describes a class of shifts which arise
from product states on $\uhf_n$, and they will be studied
and characterized further in Sections 6--7, with view to
their harmonic analysis. Our condition \rom{(5.16)} below
for {\it conjugacy\/} is closely related to \cite{Sta;
Theorem 3.6} and \cite{Lac1; Theorem 4.3}; and we are
grateful to M.~Laca for bringing the reference \cite{Sta} to
our attention.

\proclaim{Theorem 5.8} Let $(U_p)$ be a sequence of
unitaries on $\Bbb C^n$ satisfying
$$
\sum_{p=0}^\infty \|e_0 - U_p e_0 \|^2 <\infty \tag 5.9
$$
and let $T_i=S_i \Gamma (U)$ where $S_i$ is defined by
\rom{(7.12)} and $\Gamma(U) = \bigotimes_{p=0}^\infty U_p$
by \rom{(7.17)}. We have
$$
e_0 \otimes e_0 \otimes e_0 \otimes \cdots = \1\in \Cal L^2
(\Omega,\mu).\tag 5.10
$$
The state $\omega_U$ defined on $M_{n^\infty} =
\bigotimes_{m=0}^\infty M_n$ by
$$\multline
\omega_U (e_{i_1j_1} \otimes e_{i_2j_2} \otimes \cdots
\otimes e_{i_mj_m} \otimes 1 \otimes 1 \otimes \cdots)\\
=\langle \1,T_{i_1} T_{i_2} \cdots T_{i_m}T_{j_m}^*\cdots
T_{j_1}^*\1\rangle,
\endmultline\tag 5.11$$
where $e_{ij}$, $i,j=1.\ldots,n$ is a set of matrix units
for $M_n=B(\Bbb C^n)$, is a product state
$$
\omega_U = \bigotimes_{m=0}^\infty \omega_{U,m} \tag 5.12
$$
where
$$\align
\omega_{U,m}& = \langle \xi_m, \cdot \xi_m \rangle, \tag
5.13\\
\xi_0 &= e_0 \tag 5.14\\
\xi_m &=U_0^* \cdots U_{m-1}^* e_0 \tag 5.15
\endalign$$
for $m=1,2,\ldots$. Hence, if $(V_p)$ is another sequence of
unitaries on $\Bbb C^n$ satisfying conditions \rom{(5.9)--
(5.10)}, then the shift associated to $S_i\Gamma(V)$ is
conjugate to the shift associated to $S_i \Gamma(U)$ if and
only if there is a unitary $W\in U(n)$ such that
$$
\sum_{m=0}^\infty (1 - |\langle V_0^*V_1^* \cdots V_m^* e_0,
  WU_0^*U_1^* \cdots U_m^*e_0 \rangle |) <+\infty. \tag 5.16
$$\endproclaim

\demo{Proof} Note first that the condition \rom{(5.9)} is
equivalent to $\bigotimes_{p=0}^\infty U_p e_0$ being a
well-defined vector in $\bigotimes_{p=0}^\infty \Bbb C^n
=\Cal L^2 (\Omega,\mu)$ (and \rom{(5.9)} is implied by the
condition $\sum_p \| 1-U_p\|<\infty$ considered in Section
7). Thus
$$
\Gamma(U) = \bigotimes_{p=0}^\infty U_p \tag 5.17
$$
is a well-defined unitary operator on $\Cal
L^2(\Omega,\mu)$, so $T_i = S_i \Gamma(U)$ are well defined,
and
$$
\alpha \circ \gamma_U(A) =\sum_{i=1}^n T_i AT_i^* \tag 5.18
$$
for all $A\in \bh$, where $\alpha$ and $\gamma_U$ are
defined in Theorem 7.3.

It follows from \rom{(8.5)} and \rom{(8.6)} that
$$\align
S_i^* (\eta_0 \otimes \eta_1 \otimes \eta_2 \otimes \cdots)
& = n^{-1/2} \eta_0 (i) (\eta_1 \otimes \eta_2 \otimes
\cdots) \tag 5.19\\
S_i (\eta_0 \otimes \eta_1 \otimes \eta_2 \otimes \cdots)
 &= n^{1/2} (\delta_i \otimes \eta_0 \otimes \eta_2 \otimes
\cdots)\tag 5.20
\endalign$$
whenever $\eta_0 \otimes \eta_1 \otimes \eta_2 \otimes
\cdots \in \Cal L^2 (\Omega,\mu)$. Using $T_i = S_i \Gamma
(U)$, $T_i^* = \Gamma (U)^* S_i^*$, one then computes
$$\align
T_i^* (\eta_0 \otimes \eta_1 \otimes \cdots)
& = n^{-1/2} \eta_0 (i) (U_0^* \eta_1 \otimes U_1^* \eta_2
\otimes U_2^* \eta_3 \otimes \cdots)
\tag 5.21\\
T_i (\eta_0 \otimes \eta_1 \otimes \cdots)
& = n^{1/2}(\delta_i \otimes U_0 \eta_0 \otimes U_1\eta_1
\otimes \cdots).
\tag 5.22
\endalign$$
Iterating the formula for $T_i^*$, one computes
$$\multline
T_{j_m}^* T_{j_{m-1}}^* \cdots T_{j_1}^* (e_0 \otimes e_0
\otimes e_0 \otimes \cdots) \\
\shoveleft{=n^{-m/2} (U_0^* e_0)(j_2) (U_0^* U_1^* e_0)
(j_3)}\\
 \cdots (U_0^* U_1^* \cdots U_{m-2}^* e_0) (j_m)(U_0^* U_1^*
\cdots U_{m-1}^* e_0 \otimes
  U_1^* U_2^* \cdots U_m^* e_0 \otimes \cdots)
\endmultline\tag 5.23$$
and hence
$$\align
\langle \1, T_{i_1} \cdots T_{i_m}T_{j_m}^* \cdots T_{j_1}^*
\1\rangle \tag 5.24\\
& =n^{-m} \overline{U_0^* e_0 (i_2)} \cdots
     \overline{(U_0^* U_1^* \cdots U_{m-2}^* e_0)(i_m)}\\
&\qquad \cdots (U_0^* U_1^* \cdots U_{m-2}^* e_0)(j_m)
\cdots (U_0^* e_0)(j_2)\\
 &= n^{-m} \overline{\xi_0(i_1)}\, \overline{\xi_1 (i_2)}
\cdots \overline{\xi_{m-1}(i_m)}\\
 &\qquad \cdots \xi_{m-1}(j_m) \cdots \xi_0 (j_1)\\
 & = \omega_{U,0} (e_{i_1 j_1}) \omega _{U,1} (e_{i_2 j_2})
         \cdots \omega_{U,m-1} (e_{i_m j_m})
\endalign$$
if the $\xi$'s and $\omega$'s are defined by \rom{(5.12)--
(5.15)}.

Now, if we can show that $\1$ is a cyclic vector for the
representation of $\uhf_n$ defined by the $T$'s, the final
conclusion of Theorem 7.4 follows from \rom{(5.7)}. But
formuli \rom{(5.21)--(5.23)} imply
$$\multline
T_{i_1} T_{i_2} \cdots T_{i_m} T_{j_m}^* \cdots T_{j_1}^*
  (e_0 \otimes e_0 \otimes e_0 \otimes \cdots)\\
\shoveleft{ = (U_0^*e_0)(j_2)(U_0^* U_1^* e_0) (j_2)
   \cdots (U_0^* U_1^* \cdots U_{m-2}^* e_0) (j_m)}\\
  \cdots \delta_{i_1} \otimes U_0 \delta_{i_2}\otimes U_1
U_0 \delta_{i_3} \otimes
 \cdots  \otimes U_{m-2} U_{m-3} \cdots U_0 \delta_{i_m}
\otimes e_0 \otimes e_0
  \otimes e_0 \otimes \cdots
\endmultline\tag 5.25$$
The linear combinations of these vectors for a fixed $m$ are
equal to the linear combinations of vectors of the form
$\eta_0 \otimes \eta_1 \otimes \cdots \otimes \eta_{m-1}
\otimes e_0 \otimes e_0  \otimes \cdots$ and hence $\1$ is
cyclic for $T$. This ends the proof of Theorem
5.8.\qed\enddemo

\head 6. Construction of Shifts on $\bh$ With No Invariant
States \endhead
We now consider a special case, and give an explicit
construction of a family of shifts on $\bh$ which have no
pure normal invariant states. This family may be constructed
using the GNS representation theory corresponding to product
states on  $\uhf$ algebras of type $n^\infty$. A shift
$\alpha$ constructed in this manner will have Powers index
$n$, i.e., $\alpha\in \shift_n(\bh)$. This family of shifts
was already considered in Example 5.5.

We begin by fixing an integer $n\geq 2$, and then we view
$M_n(\Bbb C)$ as the algebra of linear transformations on
the $n$-dimensional vector space $\Bbb C^n$ over $\Bbb C$.
For each $k\in \Bbb N$, let $B_k$ be an isomorphic copy of
$M_n(\Bbb C)$, and in the usual way, we consider $B_k$ to be
embedded in the tensor product construction of the $\uhf$
algebra $\Cal A$ of Glimm type $n^\infty$, i.e., $\Cal
A=\bigotimes B_k$. We now construct a family of pure product
states on $\Cal A$ as follows. For each positive integer
$k$, pick a unit vector $h_k$ in $\Bbb C^n$, and let $e_k
\in M_n(\Bbb C)$ denote the corresponding rank one
projection onto $\Bbb C h_k$. Throughout this section we
impose the following conditions on the sequence of vectors
$\{h_k\}$:
$$\align
\lim_{k\rightarrow \infty}\| h_k -h \| &=0\qquad\text{for
some $h\in \Bbb C^n$}, \tag 6.1 \\
\sum_{k=1}^\infty \| h_k - h_{k+1} \| &< \infty, \tag 6.2\\
\prod_{k=m}^\infty | \langle h_k, h\rangle | &=0,\qquad
\text{for all $m$}. \tag 6.3
\endalign$$
In fact, only {\it the last two\/} are really conditions, as
\rom{(6.2)} implies that $(h_k)$ is Cauchy, and therefore
\rom{(6.1)} may be viewed as the definition of $h$. Using
the first condition, there is an integer $m$ such that
$\langle h_k,h\rangle \neq 0$ for all $k \geq m$, so the
third condition is equivalent to the divergence of the
series $\sum_{k=m}^\infty - \ln \cos |\langle h_k,h \rangle
|$. But one easily verifies that for sufficiently small $x$
in $\Bbb R$, $x^2/2 \leq -\ln \cos(x)\leq x^2$, so the
divergence of the series (and hence condition \rom{(6.3)})
is equivalent to the following condition.
$$
\sum_{k=1}^\infty \{\arccos (|\langle h_k,h\rangle |)\}^2 =
\infty. \tag6.3$'$
$$

\example{Example} We provide an example of a sequence of
vectors in $\Bbb C^2$ which satisfies \rom{(6.1), (6.2),
(6.3$'$)}. Consider the sequence of real numbers
$$\{1,1/2,1/2,1/3,1/3,1/3,\ldots\} =\{\theta_k\},$$
i.e., each term $1/q$ appears $q$ times in the sequence.
Define $h_k$ to be the vector
$[\cos(\theta_k),\sin(\theta_k)]$. Then $\{h_k\}$ converges
in norm to $h=[1,0]$, so \rom{(6.1)} is clearly satisfied.
(An alternative example would be, $\theta_k =k^{-1/2}$, for
$\forall k \in \Bbb N$.) To see that \rom{(6.2)} holds, note
that
$$\gather
\sum_{k=1}^\infty \|h_k -
h_{k+1}\|\qquad\qquad\qquad\qquad\qquad\qquad\qquad
\qquad\qquad\qquad\qquad\qquad\qquad\qquad\\
 {\align &= \sum_{q=1}^\infty \|[\cos(q^{-1}),\sin(q^{-1})]
-
        [\cos((q+1)^{-1}),\sin((q+1)^{-1})]\| \\
 &=\sum_{q=1}^\infty \{(\cos (q^{-1})- \cos((q+1)^{-1})^2
   +(\sin(q^{-1})-\sin((q+1)^{-1}))^2\}^{1/2}\\
 &= \sum_{q=1}^\infty \{2-2(\cos(q^{-1}-(q+1)^{-
1}))\}^{1/2}\\
 &=\sum_{q=1}^\infty 2^{1/2} \{1 - \cos ((q(q+1))^{-
1})\}^{1/2}
    \leq \sum_{q=1}^\infty 2^{1/2}\{1 -\cos^2 ((q(q+1))^{-
1})\}^{1/2}\\
 &=\sum_{q=1}^\infty 2^{1/2} \sin((q(q+1))^{-1})
    \leq \sum_{q=1} ^\infty 2^{1/2}/(q(q+1))
=2^{1/2}<\infty.\endalign}
\endgather$$
Finally, to see that \rom{(6.3$'$)} holds, note that
$|\langle h,h_k \rangle | = \cos (\theta_k)$, so $$
\sum_{k=1}^\infty \{\arccos (| \langle h,h_k \rangle |)\}^2
  = \sum_{k=1}^\infty \theta_k^2
  = \sum_{q=1}^\infty q(1/q)^2
  = \sum _{q=1}^\infty 1/q =\infty.
$$

For each positive integer $k$, let $\rho_k$ be the pure
state on $B_k$ defined by $\rho_k (A) =\tr(e_k A)$, where
$\tr$ is the non-normalized trace on $B_k$, and let $\rho$
be the product state $\rho= \bigotimes \rho_k$ on $\Cal A$.
The state $\rho$ is a pure state on $\Cal A$, \cite{Gui3,
Corollary 2.2}, so that in the corresponding GNS
representation $(\pi_\rho,\h_\rho,\Omega_\rho)$ for $\rho$,
the weak closure $\pi_\rho (\Cal A)''$ of $\pi_\rho(\Cal A)$
is isomorphic to $\Cal B(\h_\rho)$.

Following the development by Guichardet on infinite tensor
products of Hilbert spaces, \cite{Gui3}, (cf. \cite{vNeu})
we record some important facts about $\h_\rho$ and $\Cal
B(\h_\rho)$. Consider all formal tensor products of vectors
$x_1 \otimes x_2 \otimes \cdots$, where all but finitely
many of the vectors $x_k$ agree with the unit vectors $h_k$.
Then there is a natural inner product which is defined on
finite linear combinations of such vectors, satisfying
$$
\left\langle \bigotimes x_k,\bigotimes y_k \right\rangle
=\prod_{k=1}^\infty \langle x_k,y_k\rangle .\tag 6.4
$$
Note that all but finitely many of the inner products in the
expression for the infinite product are $1$. Then $\h_\rho$
is the Hilbert space completion, via the inner product
above, of the set of finite linear combinations of vectors
$\bigotimes x_k$, \cite{Gui3, Section 1.1} (see also
\cite{vNeu, Section 3.11} and \cite{Gui1--2}). Note that the
cyclic unit vector $\Omega_\rho$ in the GNS representation
for $\rho$ is $\bigotimes h_k$.\endexample

\proclaim{Lemma 6.1} Suppose $\{y_k : k\in \Bbb N\}$ is a
sequence of unit vectors in $\Bbb C^n$ which satisfies
$\sum_{k=1}^\infty \| y_k - h_k \| <\infty$. If for each
$p\in\Bbb N$, $H_p$ is the vector in $\h_\rho$ given by $H_p
= y_1 \otimes \cdots \otimes y_p \otimes h_{p+1} \otimes
h_{p+2} \otimes\cdots$, then $\{H_p\}$ is a Cauchy
sequence.\endproclaim

\demo{Proof} For positive integers $p<q$, $\| H_q - H_p
\|_\rho \leq \sum_{k=p+1}^q \| h_k - y_k \|$. \qed\enddemo

\remark{Remark} As a  result of the lemma it makes sense to
represent the limit of such a Cauchy sequence by the symbol
$\bigotimes y_k :=\bigotimes_{k=1}^\infty y_k$ (cf.
\cite{Gui3, Proposition 1.1}).\endremark

Next we consider the algebra $C(\h_\rho)$ of compact
operators in $\Cal B(\h_\rho)$. Our presentation is implicit
in the paper of Guichardet. For each $k\in \Bbb N$, select
matrix units $e_{ij}^k$, $1\leq i,j \leq n$, for $B_k$:
i.e., for each $k$, $e_{ij}^k e_{pq}^k = \delta_{jp}
e_{iq}^k$, and $\sum_{j=1}^n e_{jj}^k =I$. We impose the
condition $e_{11}^k=e_k$ for each $k$. For each $k$, let
$\{h_{k1},\ldots,h_{kn}\}$ be an orthonormal basis for $\Bbb
C^n$ selected so that $h_{k1}=h_k$ and $h_{kj}=e_{jj}^k
h_{kj}$, $j\in \{1,\ldots, n\}$. Next let $\Cal I$ be the
set of all ordered sequences $P=\{p_1,p_2,\ldots\}$ where
$p_k\in \{1,2,\ldots,n\}$ for each $k$, and all but finitely
many of the $p_k$ are $1$. We define $\delta_{PQ}$, $P,Q\in
\Cal I$, to be $1$ if $p_k=q_k$ for all $k$,  and otherwise
$0$. We use the notation $\bigotimes h(P)$ to represent the
unit vector $\bigotimes h_{kp_k}$ in $\h_\rho$. From the
discussion above, linear combinations of the vectors
$\bigotimes h(P)$ are dense in $\h_\rho$ and furthermore,
$\langle \bigotimes h(P),\bigotimes
h(Q)\rangle=\delta_{PQ}$, $P,Q\in \Cal I$. The following
result is clear.

\proclaim{Lemma 6.2} The set $\{\bigotimes h(P) : P\in \Cal
I\}$ forms an orthonormal basis for $\h_\rho$.\endproclaim

Next for $R,S\in \Cal I$, we use the notation $E_{RS}$ to
represent the rank one operator in $\Cal B(\h_\rho)$ which
satisfies $E_{RS}\left( \bigotimes h(P)\right)=\delta_{PS}
\left( \bigotimes h(R)\right)$. It is sometimes useful to
write $E_{RS}$ as $e_{r_1 s_1}^1 \otimes e_{r_2s_2}^2
\otimes \cdots$. From the previous equation and the previous
lemma, it follows that the operators $E_{RS}$ form a
complete set of matrix units for $C(\h_\rho)$ ($=$ the
compact operators), i.e, $E_{RS}$, satisfy the identities
$$\align
E_{PQ}^* &= E_{QP} \tag 6.5\\
E_{PQ}E_{RS} &= \delta_{QR}E_{PS},\qquad P,Q,R,S \in \Cal
I,\tag 6.6
\endalign$$
and the set of finite linear combinations of the matrix
units $E_{PQ}$ is a uniformly dense subalgebra of
$C(\h_\rho)$.

We next show that there is a natural way to make sense of
the symbol $I \otimes E_{PQ}$ as a rank $n$ operator in
$C(\h_\rho)$, and then we use these operators to define a
shift on $\Cal B(\h_\rho)$ of index $n$. To begin, let
$\bigotimes h(P)$, $P\in \Cal I$, be any vector in the
orthonormal basis for $\h_\rho$, and let $m$ be a positive
integer sufficiently large so that $p_k=1$ for any $k \geq
m$. Let $v$ be any unit vector in $\Bbb C^n$. By condition
\rom{(6.2)}, $\sum_{j=1}^\infty \|h_{m+j+1} -
h_{m+j}\|<\infty$, so by Lemma 6.1, the symbol $v \otimes
h_{1p_1} \otimes h_{2p_2} \otimes h_{3p_3} \otimes \cdots$
represents a unit vector in $\h_\rho$. Hence the symbol
$I\otimes E_{PQ}$ represents a rank $n$ operator in
$C(\h_\rho)$ which maps, for any $v\in \Bbb C^n$, the vector
$v\otimes h_{1q_1} \otimes h_{2q_2} \otimes h_{3q_3} \otimes
\cdots$ to the vector $v\otimes h_{1p_1} \otimes h_{2p_2}
\otimes h_{3p_3} \otimes \cdots$. Furthermore it is not
difficult to show that the operators $I\otimes E_{PQ}$
satisfy the identities
$$
(I\otimes E_{PQ})(I \otimes E_{RS})=\delta_{QR}(I\otimes
E_{PS}). \tag 6.7
$$
If $A\in B_1$, then clearly $\pi_\rho(A)$ and $I\otimes
E_{PQ}$ commute, for all $P,Q$ in $\Cal I$. Hence $I\otimes
E_{PQ}\in \pi_\rho(B_1)'$. On the other hand, consider the
subalgebra $\bigotimes _{k=2}^\infty B_k$ of $\Cal A$
generated by $B_2,B_3,\ldots$. From \cite{Pow1, Lemma 2.3},
$\pi_\rho(B_1)'$ coincides with $\pi_\rho(\bigotimes
B_k)''$, which is a Type I subfactor of $\Cal B(\h_\rho)$.
Thus the set of finite linear combinations of compact
operators of the form
$$
\sum_{j=1}^n e_{jj}^1 \otimes e_{r_2s_2}^2 \otimes e_{r_3
s_3}^3 \otimes\cdots
=\sum_{j=1}^n E_{R_j S_j},\tag 6.8
$$
$R,S\in \Cal I$, where for $R=\{r_1,r_2,r_3,\ldots\}$, $R_j$
is the sequence $\{j,r_2,r_3, \ldots\}$, forms a weakly
dense subalgebra of $\pi_\rho \left( \bigotimes_{k=2}^\infty
B_k \right)''$. We summarize these results below.

\proclaim{Theorem 6.3}For any $P,Q\in\Cal I$, the symbol
$I\otimes E_{PQ}$ represents a compact operator of rank $n$
in $C(\h_\rho)$. The set of such operators forms a complete
set of matrix units for the subalgebra of compact operators
of the Type I subfactor $\pi_\rho(B_1)'$ of $\Cal
B(\h_\rho)$.\endproclaim

The results of the preceding theorem enable us to define a
shift $\alpha$ of index $n$ on $\Cal B(\h_\rho)$ which
satisfies $\alpha(C(\h_\rho)) \subset C(\h_\rho)$. Fix
$S=\{1, 1,1, \ldots\}$ in $\Cal I$. Since $I\otimes E_{SS}$
is a rank $n$ projection in $C(\h_\rho)$, there exist
partial isometries $W_1,\ldots , W_n$ in $\Cal B(\h_\rho)$,
each of rank one, satisfying
$$
W_j^* W_j =E_{SS}\qquad\text{and}\qquad
\sum_{j=1}^n W_j W_j^* = I\otimes E_{SS}.
$$
Define operators $V_1,\ldots, V_n$ in $\Cal B(\h_\rho)$ by
$V_j=\sum_{K\in\Cal I} (I\otimes E_{KS})W_j(E_{SK})$ (cf.
\cite{Pow2, Theorem 2.4} it is straightforward to show that
$V_j$'s are isometries satisfying the Cuntz algebra relation
$\sum_{j=1}^n V_j^* V_j=I$. We may thus define a  shift
$\alpha$ on $\Cal B(\h_\rho)$ by setting, for $A\in \Cal
B(\h_\rho)$, $\alpha(A) = \sum_{j=1}^n V_j AV^*_j$. Note
that for $P,Q\in \Cal I$,
$$\align
\alpha(E_{PQ}) &=\sum_{j=1}^n V_j E_{PQ}V_j^*\\
  &=\sum_{j=1}^n \left\{\sum_K (I\otimes E_{KS}) W_j
E_{SK}\right\} E_{PQ}
      \left\{\sum_L E_{LS}W_j^*(I\otimes E_{SL})\right\}\\
  &=\sum_{j=1}^n (I\otimes E_{PS}) W_j^* E_{SS}W_j (I\otimes
E_{SQ})\\
  &=\sum_{j=1} ^n(I\otimes E_{PS})W_j^* W_j (I\otimes
E_{SQ})\\
  &= (I\otimes E_{PS})(I\otimes E_{SS}) (I\otimes E_{SQ})
  = I \otimes E_{PQ},
\endalign$$
so that it is natural to use the notation
$\alpha(A)=I\otimes A$ to denote this shift on $\Cal
B(\h_\rho)$. By Theorem 6.3 and the previous calculation,
$\alpha(\Cal B(\h_\rho))'= \pi_\rho(B_1)''$, so $\alpha\in
\End_n(\Cal B(\h_\rho))$, i.e., $[\Cal B(\h_\rho): \alpha
(\Cal B(\h_\rho))]=n^2$.

The following theorem gives a concrete realization of the
representation of $\on$ defined in Theorem 3.1 and Lemma 5.2
in the present setting.

\proclaim{Theorem 6.4} Let $V_1,\ldots,V_n$ be the
isometries defined as above. Then the mapping $\alpha(A)
=\sum_{j=1}^n V_j AV_j^*$ is a shift endomorphism on $\Cal
B(\h_\rho)$ of index $n$ satisfying the identities
$\alpha(E_{PQ})=I\otimes E_{PQ}$, $P,Q\in\Cal
I$.\endproclaim

\demo{Proof}To finish the proof we must show that $\alpha$
is a shift. But it is not difficult to show that
$\alpha^k(\Cal B(\h_\rho))' \supset \pi_\rho (B_1\otimes
\cdots \otimes B_k)$; and since $\left[ \bigcup_k \pi_\rho
(B_1\otimes \cdots \otimes B_k)\right]''=\Cal B(\h_\rho)$,
it follows that $\bigcap_k \alpha^k (\Cal B(\h_\rho))'=\Bbb
CI$.\qed\enddemo

Next we prove some preliminary results to be used in showing
that there are no normal $\alpha$-invariant pure states on
$\Cal B(\h_\rho)$. We could of course just refer to Lemma
5.4 and \rom{(5.7)} for this, but we prefer to give an
interesting direct argument. We state the following well-
known result for convenience.

\proclaim{Lemma 6.5} Let $H,H'$ be unit vectors in
$\h_\rho$, and let $\omega,\omega'$ be the corresponding
(pure normal) vector states on $\Cal B(\h_\rho)$. Then
$\|\omega-\omega'\|\leq 2\|H-H'\|$.\endproclaim

\demo{Proof} The result follows from the inequality
$$\align
|\omega(A)-\omega'(A)| &= |\langle AH, H\rangle - \langle
AH',H' \rangle |\\
  & \leq |\langle AH,H \rangle - \langle AH,H' \rangle | +|
\langle AH,H'\rangle - \langle AH',H' \rangle| \\
  & \leq \|AH\| \cdot \| H-H'\| + \|A(H-H')\| \cdot \|H\|.
\endalign$$\qed\enddemo

For the following two results the notation $\Cal I_m$,
$m\in \Bbb N$, is used to denote the set of sequences $P\in
\Cal I$ whose entries are all $1$ with the possible
exception of the first $m-1$ entries. Observe that $\Cal I_1
\subset \Cal I_2 \subset \cdots$ and $\bigcup_m \Cal I_m =
\Cal I$.

\proclaim{Lemma 6.6} Let $E$ be an orthogonal rank one
projection in $\Cal B(\h_\rho)$. Then for any $\epsilon>0$,
there is a finite linear combination $E'$ of the rank one
operators $E_{PQ}$, $P,Q\in \Cal I$, which is an orthogonal
projection satisfying $\|E - E' \|<\epsilon$.\endproclaim

\demo{Proof} Since the $E_{PQ}$'s, $P,Q\in \Cal I$, form a
full set of matrix units for $C(\h_\rho)$, there are, for
some $m$, sequences $P(1),P(2),\ldots,P(m),Q(1),\ldots,Q(m)$
of $\Cal I_m$, and complex numbers $c_j$, $j=1,\ldots,m$,
such that $\left\|E-\sum_{j=1}^m c_j
E_{P(j)Q(j)}\right\|<\epsilon$. Hence the sum $\sum_{j=1}^m
c_j E_{P(j)Q(j)}$ takes the form $A\otimes e_m \otimes
e_{m+1}\otimes \cdots$, for some $A\in\pi_\rho(B_1\otimes
\cdots \otimes B_{m-1})$. Using standard functional calculus
techniques \cite{Gli, Lemma 1.6} we may assume that $A$ is a
projection in $\pi_\rho(B_1 \otimes \cdots\otimes  B_{m-1})$
and the result \hbox{follows.\qed}\enddemo

\remark{Remark} Note that if $\epsilon<1$ then the
projection constructed in the proof of the previous lemma
must also be rank one.\endremark

The following lemma identifies an important {\it clustering
property\/} which we take up again in Sections 7--8 below.

\proclaim{Lemma 6.7} Let $E'$ be any projection of the form
$A\otimes e_m \otimes e_{m+1} \otimes \cdots$ as in the
previous lemma. Let $H' \in \h_\rho$ be any unit vector
obtained as a finite linear combination of vectors in the
orthonormal basis $\{h(P): P\in \Cal I\}$ of $\h_\rho$. Let
$\omega'$ be the vector state corresponding to $H'$. Then
$\lim_{k\rightarrow \infty} \omega'(\alpha^k (E'))
=0$.\endproclaim

\demo{Proof} Suppose for some $m\in \Bbb N$, $H'$ is a
finite linear combination of the vectors $\bigotimes
h(R(1)),\ldots ,\bigotimes h(R(m))$, for some sequences
$R(j)$ in $\Cal I_m$. Then we may write $H'$ in the form
$\Phi\otimes h_m \otimes h_{m+1} \otimes h_{m+2} \otimes
\cdots$, where $\Phi$ is a unit vector in the $(m-1)$-fold
tensor product of $\Bbb C^n$. From the preceding theorem,
$\alpha^k (E')$ has the form $I\otimes I\otimes \cdots
\otimes I\otimes A\otimes e_m\otimes e_{m+1} \otimes\cdots$,
where the first $k$ tensors are $I$. Note that from the form
of $H'$ and of $E'$ we have
$$\align
\omega'(\alpha^k(E'))
  & \leq \prod_{j=1}^\infty \langle
e_{m+j}h_{m+j+k},h_{m+j+k}\rangle\\
  & = \prod_{j=1}^\infty| \langle h_{m+j},h_{m+j+k}
\rangle|^2
\endalign$$
By condition \rom{(6.1)}, $\lim_{k\rightarrow
\infty}h_{m+j+k}=h$, and by condition \rom{(6.3)}, $$
\prod_{j=1}^\infty |\langle h_{m+j},h\rangle|=0.
$$
Applying these two conditions to the infinite product above,
it is not hard to show that
$$
\lim_{k\rightarrow \infty}\left\{ \prod_{j=1}^\infty|
\langle h_{m+j},h_{m+j+k}\rangle |^2 \right\} =0,
\tag 6.9
$$
and the result follows.\qed\enddemo

\proclaim{Theorem 6.8} Let $\alpha$ be a shift on $\Cal
B(\h_\rho)$ constructed as above. Then there are no pure
normal $\alpha$-invariant states on $\Cal
B(\h_\rho)$.\endproclaim

\demo{Proof} Any pure normal state $\omega$ on $\Cal
B(\h_\rho)$ is a vector state $\omega=\langle H, \cdot
H\rangle$, for some unit vector in $\h_\rho$. Given
$\epsilon>0$, there is a vector $H'$ such that $H'$ is a
finite linear combination of the basis vectors $\bigotimes
h(P)$, $P\in\Cal I$, with $\|H-H'\| <\epsilon/3$. Let $E$ be
an orthogonal rank one projection in $C(\h_\rho)$, then by
Lemma 6.6 there is a rank one projection $E'$ which is a
finite linear combination of the matrix units $E_{PQ}$,
$P,Q\in \Cal I$, such that $\|E-E'\|<\epsilon/3$. Let
$\omega'=\langle H', \cdot H'\rangle$, then $\|\omega-
\omega'\| <2\epsilon/3$, by Lemma 6.5. Then since
$\|\alpha^k(E)-\alpha^k (E')\| = \|E-E'\|$ we have, for all
$k$, $|\omega(\alpha^k(E))-\omega'
(\alpha^k(E'))|<\epsilon$. But $\lim_{k\rightarrow \infty}
\omega' (\alpha^k (E))=0$, by the previous lemma. Since
$\epsilon$ is arbitrary we have $\lim_{k\rightarrow\infty}
\omega(\alpha^k (E))=0$ also. Hence if $\omega$ were an
$\alpha$-invariant state then $\omega(E) =0$ for all
orthogonal rank one projections in $C(\h_\rho)$. But then
$\left.\omega\right|_{C(\h_\rho)}=0$  which contradicts the
normality of $\omega$. This contradiction yields the
result.\qed\enddemo

\proclaim{Corollary 6.9} Let $\alpha$ be a shift on $\Cal
B(\h_\rho)$ constructed as above. Then there are no normal
$\alpha$-invariant states on $\Cal B(\h_\rho)$.\endproclaim

\demo{Proof} Suppose $\omega$ is a finite linear combination
$\sum_{j=1}^m a_j \omega_j$, $a_j\in \Bbb R^+$, of normal
states $\omega_j= \langle H_j, \cdot H_j \rangle$. For any
$\epsilon>0$ one may, as in the proof of the theorem, choose
unit vectors $H'_j$, each of which is a finite linear
combination of the basis vectors $\bigotimes h(P)$, $P\in
\Cal I$, and satisfying $\| H_j - H_j '\| <\epsilon/3$. Then
if $\omega'=\sum_{j=1}^m  a_j  \omega'_j $, $\|\omega -
\omega'\|< 2\epsilon/3$, and for $E'$ chosen as above it is
clear that $\lim_{k\rightarrow \infty} \omega' (\alpha^k
(E'))=0$, so that $\lim_k \sup\{ |\omega (\alpha^k (E))
|\}<\epsilon$, for all rank one projections $E$ in
$C(\h_\rho)$. Since $\epsilon $ is arbitrary we therefore
obtain $\lim_{k\rightarrow \infty} \omega (\alpha^k(E))=0$,
whence, as in the proof of the theorem, $\omega$ cannot be
$\alpha$-invariant. Finally, since any normal state $\omega$
of $\Cal B(\h_\rho)$ may be approximated arbitrarily closely
in norm by states which are finite linear combinations of
vector states, we have $\lim_{k\rightarrow\infty} \omega
(\alpha^k(E)) =0$ for these states as well, so such a state
cannot be $\alpha$-invariant.\qed\enddemo

\head 7. Clustering Properties \endhead
Let $n\in\Bbb N$, $n>1$, be given, and let $\on$ be the
corresponding {\it Cuntz-algebra\/} on generators
$(s_i)_{i=0}^{n-1}$ and relations, $s_i^* s_j =\delta_{ij}
\1$, and $\sum_{i=0}^{n-1} s_is_i^* =\1$. Let $\h$ be a
separable (infinite-dimensional) complex Hilbert space. Then
we saw that each element in $\rep (\on,\h)$ is specified by
an assignment, $s_i \mapsto S_i$ of isometries of $\h$,
subject to the Cuntz-relations,
$$
S_i^* S_j =\delta_{ij} I\qquad\text{and}\qquad
\sum_{i=0}^{n-1} S_i S_i^* =I \tag 7.1
$$
where $I$ denotes the identity operator on $\h$. In Theorem
3.1, we identified the {\it $U(n)$-equivalence\/} (denoted
$\sim$) on $\rep(\on,\h)$, and a (bijective) isomorphism
$$
\End_n(\bh)\approxeq \rep(\on,\h)/ \sim. \tag 7.2
$$
The element $\alpha \in \End_n(\bh)$ which corresponds to a
given $(S_i)\in \rep(\on,\h)$ is
$$
\alpha(A) =\sum_{i=0}^{n-1} S_i A S_i^* \tag 7.3
$$
defined for $\forall A\in \bh$. We also saw in Section 2
that $\alpha$ in \rom{(7.3)} is a {\it shift\/} precisely
when the operators
$$
S_{i_1} \cdots S_{i_p} S_{j_p}^* \cdots S_{j_1}^* \tag 7.4
$$
act irreducibly on $\h$. (Note that the family in
\rom{(7.4)} is {\it indexed by\/} (variable) $p \in \Bbb N$,
and double-multi-indices, $i_1,\ldots,i_p,j_1,\ldots,j_p$.)

For any two elements $(S_i)$ and $(T_j)$ in $\rep(\on,\h)$,
it is clear from \rom{(7.1)} that the matrix
$$
(S_i^* T_j) \in \Cal M_n (\bh) \tag 7.5
$$
is {\it unitary}. Note that the matrix entries,
$M_{ij}=S_i^* T_j$ are generally just in $\bh$. It also
follows (as noted in \rom{(2.4)--(2.5)} above) that, {\it
conversely}, if $(S_i) \in \rep(\on,\h)$, and $(M_{ij})\in
\Cal M_n (\bh)$ is given unitary, {\it then\/} the operators
$T_j$ defined by
$$
T_j =\sum_i S_i M_{ij} \tag 7.6
$$
satisfy the Cuntz-relations \rom{(7.1)} and
$$
S_i^* T_j =M_{ij}. \tag 7.7
$$
We think of the unitary operator-valued matrix $(M_{ij})$ as
a {\it non-commutative Radon-Nikodym derivative\/} relating
two elements in $\rep(\on,\h)$. By \rom{(7.2)}, it will
therefore also be relating the corresponding elements in
$\End_n(\bh)$.

We will show that there is a distinguished (up to unitary
equivalence) element $(S_i)\in \rep(\on,\h)$ corresponding
to a certain Haar measure (details below). It will be  a
shift, and we shall refer to it as the {\it Haar-shift}. It
has a pure invariant state which is defined directly in
terms of the constant function on $\Omega$ (where $\Omega$
is the infinite product group defined from $\Bbb Z_n$, see
\rom{(7.11)} below) and the Haar measure on this compact
group $\Omega$ (see \rom{(7.11)} below). Our purpose in the
present section is to be able to read off from the Radon-
Nikodym derivative \rom{(7.7)} when some second element
$(T_j)\in\rep(\on,\h)$ also has a pure invariant state.
Recall, by \rom{(7.2--7.3)}, that
$$
\beta(A):=\sum_j T_j A T_j^*\qquad\text{(for $A\in \bh$)}
\tag 7.8
$$
is the element in $\End_n(\bh)$ which corresponds to the
given $(T_j)$; and that the possible existence of pure and
invariant states refers then to the possible existence of
unit-vectors $\xi\in\h$ such that
$$
\langle \xi, A\xi \rangle = \langle \xi,\beta(A)\xi \rangle
\qquad\text{for $\forall A\in\bh$.} \tag 7.9
$$
We saw in Theorem 4.1 that such a vector $\xi$ exists if and
only if there is a solution $c=(c_i)\in \ell^2_n$, with
$\sum_{i=0}^{n-1} |c_i |^2 =1$, to the {\it simultaneous
eigenvalue problem},
$$
T_j^* \xi=c_j \xi\qquad \text{for $0\leq j <n$}. \tag 7.10
$$

\definition{Definition 7.1} Following \cite{Jo-Pe}, we now
describe {\it the Haar-shift of index\/} $n$. We recall the
residue group, $\Bbb Z_n := \Bbb Z/n \Bbb Z \simeq\{
0,1,\ldots, n-1\}$, and the corresponding infinite Cartesian
product group,
$$
\Omega=(\Bbb Z_n)^{\Bbb N} = \prod_{p=1}^\infty \Bbb Z_n.
\tag 7.11
$$
It is viewed as a compact abelian group under coordinate
addition. The corresponding normalized Haar measure on
$\Omega=\Omega_n$ will be denoted $\mu$. It is the product
measure corresponding to assigning equal weights $n^{-1}$ at
the $n$ coordinates (of each factor). Points in $\Omega$ are
denoted
$$
x=(x_p)_{p=1}^\infty =(x_1,x_2,\ldots),
$$
and we have the right and left Bernoulli-shifts given
respectively by
$$
\sigma_i(x_1,x_2,\ldots)  =(i,x_1,x_2,\ldots),
\qquad
\text{and}
\qquad
\sigma(x_1,x_2,\ldots)=(x_2,x_3,\ldots).
$$
Clearly then, $\sigma\circ \sigma_i =\id_\Omega$ for all
$i$, and furthermore,
$$
\mu=\frac 1n \sum_{i=0}^{n-1} \mu \circ \sigma_{i}^{-
1},\qquad \mu\circ\sigma_i=\frac 1n \mu,
$$
and therefore $\mu\circ \sigma^{-1} =\mu$.

It follows (see \cite{Jo-Pe}) that we get a Cuntz-algebra
system $(S_i)$ on $\h=\Cal L^2 (\Omega,\mu)$ as follows: The
operators $S_i$, and their adjoints $S_i^*$, will be acting
on $\h$, and are given by,
$$
S_i^* \xi =n^{-1/2} \xi \circ \sigma_i \qquad\text{for
$\forall \xi\in\h=\Cal L^2(\Omega,\mu)$.}
\tag 7.12
$$
The corresponding shift $\alpha$ from \rom{(7.3)} will be
called {\it the Haar-shift}. The vector state $\omega_0$ on
$\bh$, given by Haar-measure $\mu$, and the constant
function $\1$, is $\alpha$-invariant: For $A\in \bh$, we
therefore have,
$$\align
\omega_0 (A) &=\langle \1,A\1 \rangle_{\Cal L^2}
=\int_\Omega (A\1)(x)\, d\mu(x), \tag 7.13a\\
\omega_0(\alpha(A))&=\omega_0(A). \tag 7.13b
\endalign$$

We shall need the {\it character\/} on the group $\Bbb Z_n$,
defined as follows: For  $p\in \Bbb Z\,$, set
$$
e(p) := \exp(i 2\pi p/n), \tag 7.14
$$
and, for $j\in \Bbb Z_n$, $k\in \Bbb Z_n$, $e(jk)$ is given
by this, with $p=jk$ and $jk$ representing the product in
the ring $\Bbb Z_n$. We shall write, $e_j(k) := \langle j,k
\rangle = e(jk)$.\enddefinition

\definition{Definition 7.2} To be able to describe our
Radon-Nikodym derivative, we shall need a certain {\it
unitary representation\/} acting on $\h=\Cal
L^2(\Omega,\mu)$.

Consider first the infinite product of identical copies of
the group $U(n)$ of all $n$ by $n$ unitary matrices. Inside
this product, we have the infinite-dimensional {\it
subgroup\/} of elements $U=(U_p)_{p=1}^\infty$ subject to,
$$
\sum_{p=1}^\infty \|I-U_p\| <\infty \tag 7.15
$$
where $I$ is the $n\times n$ identity matrix $\pmatrix
1&0&\cdots &0\\
0&1& & \\
\vdots& &\ddots & \\
0& & &1\endpmatrix$, and $\|\cdot\|$ is the $C^*$-norm on
the $n\times n$ matrices. (In fact the weaker condition
$\sum_p \| I-U_p \|^2 <\infty$ will suffice.) This subgroup
will be denoted $G_n$, and it has a natural {\it unitary
representation\/} on $\Cal L^2 (\Omega,\mu)$ which we
proceed to describe.

Using \rom{(7.14)}, we note that the discrete group
$\Lambda$, which is {\it dual\/} to $\Omega = \prod \Bbb
Z_n$, is $\Lambda = \coprod \Bbb Z_n$ consisting of
elements, $\lambda= (y_1, \ldots, y_q, 0,0,\ldots)$ with at
most a finite number of nonzero points $y_j$ in $\Bbb Z_n$,
followed by an {\it infinite\/} string of zeros. We get an
{\it orthonormal basis\/}  $e_\lambda$, indexed by
$\lambda\in\Lambda$, and given by,
$$
e_\lambda (x) := \prod _{p=1}^q e(y_p x_p)
=\prod_{p=1}^q \langle y_p, x_p \rangle. \tag 7.16
$$Note that we may also view \rom{(7.16)} as an {\it
infinite product}, but the factors after $q$ will all be
one. For $y\in \Bbb Z_n$, we further have the functions $e_y
\in \ell^2_n = \ell^2 (\Bbb Z_n) \simeq \Bbb C^n$, given by,
$e_y(x):= e(yx)$, see \rom{(7.14)} above. This is again an
orthonormal basis, now in $\ell_n^2$, relative to the Haar
measure on $\{0,\ldots,n-1\}$, i.e., equal weights $n^{-1}$.
Each unitary $n$ by $n$ matrix may then be identified with a
unitary transformation on $\ell_n^2$ relative to this basis.
For an element, $U=(U_p)_{p=1}^\infty$ in $G$, we define
$\Gamma(U)$ on the basis $\{e_\lambda\}_{\lambda\in\Lambda}$
as follows:
$$
(\Gamma(U)e_\lambda)(x) :=
\prod_{p=1}^\infty (U_p e_{y_p})(x_p). \tag 7.17
$$
Using the argument from Section 6, we may then check that
the right-hand side in \rom{(7.17)} represents a well-
defined element in $\Cal L^2 (\Omega,\mu)$, with the
infinite product {\it convergent in mean-square}. We omit
the simple argument which is based directly on the
summability \rom{(7.15)} defining the subgroup $G$. It is
also clear that, $U \rightarrow \Gamma(U)$, is then a
unitary representation of $G$ acting on $\Cal
L^2(\Omega,\mu)$.
\enddefinition

The construction of the unitary representation, $U \mapsto
\Gamma(U)$ in \rom{(7.17)} is parallel to the corresponding
one, $U\mapsto \widetilde{\Gamma}(U)$ from \cite{Gui2--3};
but with the $\widetilde{\Gamma}$ representation acting on
von Neumann's Hilbert space $\h(h_p)$ associated with some
(fixed) sequence $(h_p)_{p=1}^\infty$ (specified as in
Section 6; see especially Lemma 6.6 and formula \rom{(6.4)}
for details): When $U=(U_p)_{p=1}^\infty$ in $G$ is given,
then $\widetilde{\Gamma}(U)$ is defined on the generic
monomials in $\h(h_p)$ by the ansatz:
$$
\widetilde{\Gamma}(U) \bigotimes_{p=1}^\infty \eta_p :=
\bigotimes_{p=1}^\infty U_p \eta_p.
$$
We shall need below a specific unitary isomorphism defined
in \rom{(7.24)}
$$
W:\h(h_p)\overset \approx \to \rightarrow \Cal
L^2(\Omega,\mu)
$$
which {\it intertwines\/} the two representations, i.e., we
have
$$
\Gamma(U) W= W\widetilde{\Gamma}(U)\qquad\text{for $\forall
U\in G$}.
$$

We are now ready to state the main result of the present
section.

\proclaim{Theorem 7.3} Let $n\in \Bbb N$, $n>1$, be given,
and let $\Omega$ be the corresponding infinite product
\rom{(7.11)}. Let $\Cal H=\Cal L^2(\Omega,\mu)$, and let
$$
\alpha(A) = \sum_i S_i A S_i^*,\qquad A\in \bh
$$
be the Haar-shift from \rom{(7.12)}. Let
$U=(U_p)_{p=1}^\infty \in G_n$ be given such that, for all
$a,b\in \ell_n^2$ of unit-norm, i.e., $\sum |a_i|^2 =\sum
|b_i|^2 =1$, we have
$$
\sum_{p=1}^\infty \left( \cos^{-1} |\langle a, U_1 \cdots
U_p b\rangle |\right)^2=\infty. \tag 7.18
$$
Let $\gamma =\gamma_U \in \aut(\bh)$ be given by
$$
\gamma_U (A) =\Gamma (U)A\Gamma (U)^* \qquad\text{for
$\forall A\in \bh$}. \tag 7.19
$$
Then
$$
\beta:= \alpha \circ \gamma_U \tag 7.20
$$
 is a shift of multiplicity $n$ which has no pure normal
invariant states. Moreover, we have
$$
\beta(A)= \sum_j T_j A T_j^*\qquad\text{for $\forall A\in
\bh$} \tag 7.21
$$
where
$$
T_j:= S_j\Gamma(U)\qquad\text{(for $\forall j$)}.
$$\endproclaim

\demo{Proof} By Theorem 3.3, the endomorphism $\beta$ in
\rom{(7.21)} may be defined alternatively from an element
$(T_j ')\in \rep(\on,\h)$ with corresponding Radon-Nikodym
derivative,
$$
S_j^* T_k ' = n^{-1/2} e(jk) \Gamma (U); \tag 7.22
$$
and this representation is the one we identify (up to
unitary equivalence) in Section 6 above, but acting in von
Neumann's infinite-product Hilbert space. The result then
follows from our Theorems 3.1 and 6.8 above. Let
$H:=\ell_n^2$, and let $(h_p)_{p=1}^\infty$ be a sequence of
vectors in $H$ such that $\| h_p \| =1$, and $\exists h\in
H$ such that,
\roster\item"{(i)}" $\lim_{p \rightarrow \infty} h_p =h$,
\item"{(ii)}" $\sum_{p=1}^\infty \| h_p - h_{p+1} \| <
\infty$ (recall that (i) is implied by (ii)), and
\item"{(iii)}" $\sum_p\left(\cos^{-1} |\langle h_p, h\rangle
| \right)^2 =\infty$.\endroster
Then we saw in Theorems 6.4 and 6.8  that von Neumann's
Hilbert space $\Cal H(h_p)$ (specifics in \cite{vNeu})
carries a shift $\tilde \beta(A) =I\otimes A$ which has {\it
no\/} pure invariant states. If $v_i\in H$ is an orthonormal
basis, and
$$
\tilde T_i \xi := v_i \otimes \xi \qquad\text{for $\forall
\xi \in \h(h_p)$}, \tag 7.23
$$
then $(\tilde T_i ) \in \rep(\on,\h(h_p))$ and $\tilde \beta
(A) =\sum_i \tilde T_i A \tilde T_i^*$ for $\forall A\in
\Cal B(\h(h_p))$. If $v\in H$ is chosen such that $V_p
v=h_p$ for a sequence of unitaries $(V_p)_{p=1}^\infty$,
then the unitaries, $U_p := V_p V_{p+1}^*$, satisfy, $h_p
=U_p h_{p+1}$; and we have a correspondence between our
conditions (i)--(iii) on the one hand, and the two
conditions (7.15) and (7.18) for the sequence
$(U_p)_{p=1}^\infty$ on the other hand. (Note that
\rom{(7.18)} is equivalent to $\sum_p(1- |\langle a,
U_1\cdots U_pb \rangle |)=\infty$.) But, if $(T_j ') \in
\rep(\on, \Cal L^2 (\Omega,\mu))$ is given by \rom{(7.22)},
and $(\tilde T_j) \in \rep (\on, \h(h_p))$ by \rom{(7.23)},
{\it then\/} we can show that they are intertwined by a
unitary isomorphism, $W: \h(h_p) \rightarrow \Cal L^2
(\Omega,\mu)$. To describe $W$, pick, for each $p\in \Bbb
N$, an orthonormal basis $(b_{j_p}^{(p)})$, indexed by $j_p
\in \Bbb Z_n$, such that $b_0^{(p)}=h_p$; and, using Lemma
6.2, we get an associated orthonormal basis,
$$
b(\lambda):= \bigotimes_{p=1}^q b_{j_p}^{(p)} \otimes
\bigotimes_{i=q+1}^\infty h_i \tag 7.24
$$
for $\Cal H(h_p)$. We then define our $W$ by sending the
basis element $b(\lambda)$ in \rom{(7.24)} to $e_\lambda \in
\Cal L^2 (\Omega, \mu)$, corresponding to the $\Lambda$-
index given by $\lambda = (j_1,\ldots, j_q, 0, 0, \ldots)$;
and it can easily be checked now that $W: \h(h_p) {\overset
\approx \to \rightarrow } \Cal L^2 (\Omega,\mu)$ has the
stated intertwining property, i.e., that $T_i 'W=W\tilde
T_i$ for $\forall i$. The proof is completed.\qed\enddemo

\remark{Remarks 7.4} (i) The fact that $\beta$ from
\rom{(7.21--7.22)} satisfies \rom{(7.20)} follows from
substitution of $T_k ' =n^{-1/2} \sum_j e(jk) S_j \Gamma
(U)$ into \rom{(7.21)}, (in fact also directly from Theorem
3.3 with $g=[n^{-1/2}e(jk)]_{jk}$):
$$\align
\beta(A)&=n^{-1} \sum_{j_1}  \sum_{j_2} \sum_k e(kj_1)
\overline{e(kj_2)}
  S_{j_1} \Gamma(U) A \Gamma (U)^* S_{j_2}^*\\
  &=\sum_j S_j \gamma_U (A)S_j^* = \sum_jT_jAT_j^*\\
  &=\alpha (\gamma_U(A))\qquad\text{for $\forall A\in \Cal
B(\Cal L^2(\Omega,\mu))$}.
\endalign$$

(ii) Let $n\in \Bbb N$, $n>1$, be given and let $G_n$ denote
the subgroup in $\prod_1^\infty U(n) =U(n)^{\Bbb N}$ defined
by condition \rom{(7.15)} above. Let $\alpha$ denote the
Haar-shift of $\Cal B(\Cal L^2(\Omega_n))$. Theorem 7.3 is
then the assertion that $\{ \alpha \circ \gamma_U : U\in
G_n\}$ contains more than one conjugacy class, so it makes
explicit the analysis from \cite{Pow2, Theorem 2.3}. We
showed that, when $U$ in $G_n$ satisfies \rom{(7.18)}, then
$\alpha \circ \gamma_U$ represents a conjugacy class {\it
different from\/} that of $\alpha$.

(iii) In Example 5.7 and Theorem 5.8, we gave a complete
abstract labeling of all the conjugacy classes of shifts
considered in the present section. The labeling is the set
of tensor products $\bigotimes_p U_p$, $\bigotimes_p V_p$
satisfying \rom{(5.9)} (or stronger \rom{(7.15)}) modulo the
equivalence relation $\bigotimes_p U_p \sim \bigotimes_p
V_p$ defined by \rom{(5.16)}. This labeling is non-smooth,
as we may expect from Theorem 1.1, and there is a continuum
of distinct conjugacy classes of this form. In Example 5.6,
we singled out subsets of the conjugacy classes in
$\shift_n(\bh)$ which were labeled by the points in a smooth
manifold. Otherwise, the other classes we have considered in
Examples 5.5 and 5.7 and Sections 6 and 7 (which are all the
same except for the difference between \rom{(5.5)} and
\rom{(6.1)}, and between \rom{(5.9)} and \rom{(7.15)}) do
not allow a complete smooth labeling. It would be
interesting to understand how numerical labels separating
conjugacy classes of $n$-shifts may possibly be assigned,
like the clustering labels in \rom{(7.25)} and \rom{(7.26)}
below. The situation is somewhat analogous to that in von
Neumann factors. Von Neumann had a discrete labeling (I$_n,
n=1,2,\ldots ,\infty,$ II$_1,$ II$_\infty$, and III). In
1967 Powers introduced a real label $\lambda$ to distinguish
isomorphism classes III$_\lambda$ in III, $0\leq \lambda
\leq 1$, and Connes and Takesaki introduced a non-smooth
label, the flow of weights, to distinguish III$_0$ classes.
A modest attempt of introducing some continuous labels is
done in Remark 7.6. The set $P/\sim$ from Theorem 1.1 and
Section 5 above (see especially details in Example 5.7) {\it
is\/} in fact a complete labeling of the $n$-shift conjugacy
classes. We showed also that {\it some of\/} the labels for
$n$-shift conjugacy classes may be identified as points in
our group $G_n$, but there are certainly other labels as
well. We will encounter one of them in Section 8.\endremark

To stress the difference between the two conjugacy classes
represented by the Haar-shift $\alpha$, and by $\beta_U :=
\alpha \circ \gamma_U$, from Theorem 7.3 above, we include
the following:

\proclaim{Corollary 7.5} Let $n\in \Bbb N$, $n>1$, be given.
Let $\alpha$ be the corresponding Haar-shift, and let $U\in
G_n$ be given subject to \rom{(7.18)}, and let $\beta_U :=
\alpha \circ \gamma_U$ be the corresponding transformed $n$-
shift. Then we have, for all $A\in C(\Cal L^2
(\Omega_n,\mu))$ ($=$ the compact operators) and all $\xi\in
\Cal L^2(\Omega_n,\mu)$, the two limits
$$
\lim_{k\rightarrow\infty} \langle \xi,\alpha^k (A)\xi
\rangle
=\omega_0 (A) \| \xi \|^2 \tag 7.25
$$
and
$$
\lim_{k\rightarrow \infty} \langle \xi, \beta_U^k (A)\xi
\rangle
=0. \tag 7.26
$$\endproclaim

\demo{Proof} We have already noted that \rom{(7.26)} is
contained in the proof of our Lemma 6.7 and Theorem 7.3
above. We showed that the problem for $\Cal
L^2(\Omega_n,\mu)$ was equivalent to one in the von Neumann
tensor product space $\Cal H(h_p)$ for a certain sequence
$(h_p)_{p=1}^\infty$ of vectors in $\ell_n^2$; and we
checked \rom{(7.26)} in $\h(h_p)$ in Theorem 6.8 (and
especially Lemma 6.7) by an approximation both in $\xi$ and
$A$.

Formula \rom{(7.25)}, on the other hand, may be checked
directly from \rom{(7.8)}, and an iteration of the formula
\rom{(7.12)} for $S_i^*$. We recall from \rom{(7.13)} that
$\omega_0(\cdot)$ is calculated directly from the Haar-
measure $\mu$ on $\Omega_n$. We omit further details on
\rom{(7.25)}, but refer instead to the paper \cite{Jo-
Pe}.\qed\enddemo

\remark{Remark 7.6} Formula \rom{(7.25)}, and recent ideas
from \cite{Pow3}, suggest the possibility of other conjugacy
invariants for $\shift_n(\bh)$. If $\alpha$ is a shift of
index $n$ on $\bh$, then any weak limit point of the
sequence $(\alpha^m(A))$ as $m\rightarrow \infty$ has to be
in $\bigcap_m \alpha^m (\bh) = \Bbb C 1$ for all $A\in\bh$,
and hence all weak limit points are scalar multiples of $1$.
Thus, if $\delta$ is any free ultrafilter on $\Bbb N$, we
may define a state $\omega(\delta)$ on $\bh$ by
$$
\omega(\delta)(A)1 = \underset {N\rightarrow \delta}\to{w -
\lim}\,A_N
$$
where
$$
A_N = \frac1N \sum_{m=0}^{N-1} \alpha^m (A),
$$
and then, of course,
$$
\lim_{N \rightarrow \delta}\frac 1N \sum_{m=0}^{N-1} \langle
\xi,\alpha^m(A)\xi \rangle
  = \omega(\delta)(A) \| \xi\|^2.
$$
As $\|\alpha(A_N) -A_N\| \leq \frac2N \|A\| \rightarrow 0$
for $N\rightarrow\infty$, the state $\omega(\delta)$ is then
$\alpha$-invariant.

If there is a state $\omega$ on $\bh$ such that $\langle
\xi,\alpha^m(A)\xi \rangle$ tends to $\omega(A)\|\xi\|^2$ in
any stronger sense, for example
$$
\lim_{m\rightarrow \infty} \langle \xi, \alpha^m(A)\xi
\rangle
  =\omega(A) \| \xi \|^2 \tag 7.27
$$
or
$$
\lim_{N \rightarrow\infty}\frac 1N \sum_{m=0}^{N-1} \langle
\xi, \alpha^m(A)\xi \rangle
  =\omega(A) \| \xi\|^2 \tag 7.28
$$
then $\omega(\delta)=\omega$, independently of $\delta$.

Now if $\alpha$ has an invariant vector state, then this
state is a Cuntz state in restriction to the representation
$\pi$ of $\on$ defining $\alpha$, by Theorem 4.1. If
$\omega$ denotes the normal extension of this state to
$\bh$, then
$$
\lim_{m\rightarrow\infty} \langle \xi,\alpha^m (A)\xi
\rangle
  =\omega(A) \|\xi\|^2
$$
for all $A\in C(\h)$, by the same reasoning as in Corollary
7.5. (This is the absorption property of \cite{Pow3}.)  But
any state on $C(\h)$ ($=$ the compact operators), has a
unique extension as a state to $\bh$, and hence
$\omega(\delta)=\omega$ for any $\delta$, also in this case.
On the other hand, if $\alpha$ does not have invariant
vector states, then $\omega(\delta)$ is necessarily non-
normal.

Note also that if $x\in\uhf_n$, and $\xi\in \h$ with $\|\xi
\| =1$, then
$$
\lim_{m\rightarrow \infty} |\langle \xi,\alpha^m
(\pi(x))\xi\rangle
 - \langle \xi, \alpha^{m+1}(\pi(x))\xi\rangle | =0
$$
by Lemma 5.2, and hence, if $\omega_0(\delta)$ is defined on
$\bh$ by
$$
\omega_0 (\delta) (A) I =\underset{m\rightarrow \delta}\to{w
- \lim}\, \alpha^m(A),
$$
then $\omega_0(\delta)$, restricted to the weakly dense
subalgebra $\pi(\uhf_n)$ of $\bh$, is $\alpha$-invariant,
and clearly $\omega(\delta)$ is an extension of
$\omega_0(\delta)$ from $\pi(\uhf_n)$ to $\bh$, i.e.,
$$
\omega(\delta) (\pi(x)) I =\underset{m\rightarrow
\delta}\to{ w - \lim}\, \alpha ^m (\pi(x))
$$
for $x\in \uhf_n$. If we put
$\omega(\alpha,\delta)=\omega(\delta)$, and if
$\gamma\in\aut(\bh)$, it is easily verified that
$$
\omega(\gamma \alpha \gamma^{-1},\delta)
=\omega(\alpha,\delta)\circ \gamma^{-1}.
$$
It is presently unclear how to get a conjugacy invariant out
of this, and relate this invariant to $P/\sim$. On the other
hand, we are able to verify the absorption property
\rom{(7.27)} for a class of shifts related to those
considered in the previous section. (For more on the
absorption property, see \cite{Pow3; Definition 2.12}.)  In
the following result we return to the notation of Section
6.\endremark

\proclaim{Theorem 7.7} Suppose $\{h_k: k\in \Bbb N\}$ is a
sequence of unit vectors in $\Bbb C^n$ satisfying the
conditions
\roster\item"{(i)}" $\sum_{k=1}^\infty \| h_{k+1} - h_k \| <
\infty$,
\item"{(ii)}" $\lim_{m\rightarrow \infty} \prod_{k=m}^\infty
\langle h_k, h \rangle =1$,
\endroster where $h=\lim_{k \rightarrow\infty} h_k$. Let
$\rho =\bigotimes \rho_k$ be the pure product state on
$\uhf_n$ where $\rho_k = \langle h_k,\cdot h_k \rangle$,
where GNS representation $(\pi_\rho, \h_\rho, \Omega_\rho)$.
Let $\omega$ be the symmetric pure product state $\omega =
\bigotimes \omega_h$ on $\uhf_n$, $\omega_h = \langle
h,\cdot h \rangle$. Then $\lim_{k\rightarrow \infty} \langle
\xi, \alpha^k (A)\xi \rangle = \omega(A) \| \xi \| ^2$, for
all $A \in \Cal B(\h_\rho)$ and all $\xi\in \h_\rho$, where
$\alpha$ is the shift given by $\alpha(A)=I\otimes A$ on
$\Cal B(\h_\rho)$.\endproclaim

\remark{Remark 7.8} Note that condition (i) is the same as
\rom{(6.2)} above. If one assumes that condition (i) holds,
then (ii) is the negation of condition
\rom{(6.3)}.\endremark

\demo{Proof}We first recall that only condition \rom{(6.2)}
was used in the proof of Theorem 6.4 so that there exists a
shift $\alpha$ on $\Cal B(\h_\rho)$ which satisfies
$\alpha(A) = I \otimes A$ for all $A$ in $\Cal B(\h_\rho)$.

Next, since for sufficiently large $m$, $\prod_{k=m}^\infty
\langle h_k, h \rangle$ exists and is nonzero, it follows
(\cite{vNeu, Lemma 2.4.2}) that $\prod_{k=m}^\infty |\langle
h_k, h \rangle |$ exists and also $\sum_{k=1}^\infty
|\theta_k| <\infty$, where $\theta_k \in (-\pi,\pi]$ is the
argument of $\langle h_k, h \rangle$. Hence, since
$|e^{i\theta}-1| \leq |\theta|$ for $\theta \in (-\pi,\pi]$,
$$\align
\sum_{k=1}^\infty |1- \langle h_k, h \rangle |
& \leq \sum_{k=1}^\infty \left\{ \left| \langle h_k, h
\rangle \right| \cdot |e^{i\theta_k} -1 |
   +\left| |\langle h_k, h \rangle | -1\right| \right\}\\
& \leq \sum_{k=1}^\infty \left\{ |\theta_k| + \left|
|\langle h_k,h\rangle | -1\right|\right\}\\
& < \infty,
\endalign$$
so (\cite{Gui3, Proposition 1.1}), $h\otimes h \otimes
h\cdots$ represents a unit vector in the Hilbert space
$\h_\rho$. For simplicity we write $H=h\otimes h\otimes
\cdots$.

Now suppose that $\xi$ is a unit vector in $\h_\rho$, then
arguing as in Lemma 6.7 there is a positive integer $m$ and
a unit vector $\xi'$ which is a finite linear combination of
vectors among the orthonormal set $\{\bigotimes h(P): P\in
I_m\}$, and satisfying $\|\xi - \xi' \|<\epsilon/4$. Write
$\xi' =\sum_{P\in I_m} a_P (\bigotimes h(P))$. The maximum
number of nonzero terms in this sum is $N=n^m$. Then using
the fact that the vector $H$ lies in $\h_\rho$, one  may
show that there exists a positive integer $M>m$ sufficiently
large so that if for each $P\in I_m$ one obtains a new
vector $H_P$ from $\bigotimes h(P)$ by replacing the tail
$\otimes h_{M+1} \otimes h_{M+2} \otimes\cdots$ of
$\bigotimes h(P)$ with $\otimes h\otimes h \otimes\cdots$,
then $\left\| \bigotimes h(P) -H_P \right\| <
\epsilon/(4N)$. Then if $\xi''$ is the vector $\sum_{P\in
I_m} a_P H_P$, one sees that $\xi '' $ is a unit vector
satisfying $\| \xi' - \xi'' \| \leq \sum_{P\in I_m}  |a_P|
\, \left\| \bigotimes h(P) -H_P \right\| \leq N \cdot 1\cdot
\epsilon/(4N) = \epsilon/4$. Hence $\|\xi - \xi'' \| <
\epsilon/2$, and therefore, by Lemma 6.5, $|\langle
\xi,\alpha^k (A)\xi \rangle - \langle \xi'', \alpha^k (A)
\xi'' \rangle|\leq \epsilon \| A\|$, for all $A\in \Cal
B(\h_\rho)$ and all $k\in \Bbb N$.

But if $k$ is chosen to be greater than $M$, note that
$\langle \xi'',\alpha^k(A)\xi'' \rangle = \langle H,AH
\rangle = \omega(A)$. Since $\epsilon$ is arbitrary, we
obtain $\lim_{k\rightarrow \infty} \langle \xi,
\alpha^k(A)\xi \rangle =\omega(A)$.\qed\enddemo

\head 8. Nearest Neighbor States \endhead
In Sections 6 and 7, we constructed shifts on $\bh$ coming
from product states on $\uhf_n$. In this section, we will
consider a state on $\uhf_n$ which is  a prototype of what
could be called a nearest neighbor state, since it couples
nearest neighbors in the tensor product decomposition
$\uhf_n = M_n \otimes M_n \otimes M_n \otimes \cdots$. We
will study this shift by perturbing the shifts with
invariant states considered in Section 4. To this end we
need to describe the latter more explicitly. We assume $n\in
\{2,3,\ldots\}$.

Let $\eta =(\eta_0, \eta_1, \ldots, \eta_{n-1})$ be a
sequence of complex numbers with
$$
\sum_{k=0}^{n-1}  |\eta_k |^2 =1.
$$
We also assume for the moment that $\eta_k \neq 0$ for
$k=0,\ldots,n-1$. Let $\Omega =\prod_{k=0}^\infty \Bbb Z_n$,
so that $\Omega$ is homeomorphic to the Cantor set. Equip
$\Omega$ with the infinite product measure $\mu$ obtained
from the measure on $\Bbb Z_n$ with weights $|\eta_0|^2,
|\eta_1|^2, \ldots, |\eta_{n-1}|^2$ on the $n$ points.
Define continuous open injections $\sigma_i: \Omega
\rightarrow \Omega$ by
$$
\sigma_i(x_0, x_1, x_2, \ldots)=(i, x_0, x_1,\ldots) \tag
8.1
$$
and define the shift $\sigma:\Omega\rightarrow \Omega$ by
$$
\sigma (x_0,x_1,x_2,\ldots)=(x_1,x_2,x_3,\ldots). \tag 8.2
$$
The corresponding element in $\rep(\on,\Cal L^2(\mu))$ may
be identified by: Define
$$\align
S_i^* \xi &= \bar\eta_i \xi \circ \sigma_i \tag 8.3\\
S_i \xi &= \bar\eta_i^{-1} \chi_{\sigma_i \Omega}\xi\circ
\sigma, \tag 8.4
\endalign$$
or
$$\align
(S_i^*\xi)(x_0,x_1, x_2, \ldots) &= \bar\eta_i \xi (i,x_0,
x_1,\ldots) \tag 8.5\\
(S_i\xi)(x_0,x_1,x_2,\ldots)&=\bar\eta_i^{-1}
\delta_{ix_0}\xi(x_1,x_2,\ldots). \tag 8.6
\endalign$$
One checks, using the formula (see \cite{Kak})
$$
\int_\Omega \psi(x_0, x_1, \ldots)\,d\mu (x_0, x_1, \ldots)
=\sum_{i=0}^{n-1} |\eta_i|^2 \int_\Omega \psi
(i,x_1,x_2,\ldots)\,d\mu (x_1,x_2,\ldots),
$$
that $S_i^*$ is indeed the adjoint of $S_i$, and that $S_i$
satisfy the Cuntz relations \rom{(2.1)}.

In fact notice that {\it distinct\/} weight distributions,
$p=(p_i)_{i\in \Bbb Z_{n}}$, where $p_i := |\eta_i|^2 >0$,
give corresponding {\it orthogonal\/} (i.e., mutually
singular) measures $\mu=\mu_p$ on $\Omega=\prod_0^\infty
\Bbb Z_n$ by an application of Kakutani's theorem
\cite{Kak}. However the individual operators $S_i$ in
\rom{(8.6)} depend both on the $p_i$'s and on the phases
$\eta_i |\eta_i|^{-1}$.  Note also that the constant
function $\xi=\1$ is a joint eigenvector for $S_1^*, \ldots,
S_n^*$ with eigenvalues $\bar\eta_1,\ldots, \bar\eta_n$, and
hence $\langle \1,\cdot \1\rangle$ defines the Cuntz state
on $\on$ by Theorem 4.1. We have
$$
(S_{i_1}\cdots S_{i_k} \1)(x_0,x_1,x_2,\ldots)
= \bar\eta_{i_1}^{-1} \delta_{i_{1},x_{0}} \bar\eta_{i_2}^{-
1} \delta_{i_{2},x_{1}}
  \cdots \bar\eta_{i_k}^{-1} \delta_{i_{k},x_{k-1}} \tag 8.7
$$
and hence $\1$ is a cyclic vector for the representation.

Note further that
$$\multline
(S_{i_1} \cdots S_{i_k} S^*_{j_k} \cdots S_{j_1}^*
\xi)(x_0,x_1,\ldots)\\
=\bar\eta_{i_1}^{-1} \delta_{i_{1},x_{0}} \cdots
\bar\eta_{i_k}^{-1} \delta_{i_{k},x_{k-1}}
\bar\eta_{j_k} \cdots \bar\eta_{j_1} \xi(j_1 , \ldots, j_k,
x_k,x_{k+1},\ldots)
\endmultline\tag 8.8
$$
and hence
$$\split
(S_{i_1} \cdots S_{i_k} S_{i_k}^* \cdots S_{i_1}^* \xi)
(x_0,x_1,\ldots)\\
&=\delta_{i_{1},x_{0}}\delta_{i_{2},x_{1}} \cdots
\delta_{i_{k},x_{k-1}}
  \xi (i_1,\ldots i_k,x_k,x_{k+1}, \ldots)\\
&=\delta_{i_{1},x_{0}} \delta_{i_{2},x_{1}} \cdots
\delta_{i_{k},x_{k-1}} \xi
  (x_0,x_1,x_2,\ldots).
\endsplit\tag 8.9
$$

It follows from \rom{(8.8)} that $\uhf_n$ acts irreducibly
on $\Cal L^2(\Omega,\mu)$, confirming by Theorem 3.3 that
the corresponding endomorphism of $\bh$ is a shift. It
follows from \rom{(8.9)} that $\pi(D_n)$ identifies with
$C(\Omega)$ acting as multiplication operators on $\Cal
L^2(\Omega,\mu)$.

Since, as we have seen in the proof of Theorem 4.2, the
canonical action of $U(n)$ acts transitively on the Cuntz
states, one may obtain concrete realizations of the
representation associated to a unit vector
$(\eta_0,\ldots,\eta_{n-1})$ in $\Bbb C^n$, where some of
the components are zero, by applying canonical actions on
states where all the components are nonzero.

For simplicity, let us specialize to the case
$$
\eta_i= n^{-1/2};\qquad i=0,\ldots,n-1, \tag 8.10
$$
so
$$\align
S_i^*\xi& = n^{-1/2} \xi \circ \sigma_i \tag 8.11\\
S_i\xi& = n^{1/2} \chi_{\sigma_{i}\Omega} \xi \circ \sigma
\tag 8.12
\endalign$$
for $\xi\in \Cal L^2 (\Omega,\mu)$.

In this case $\Cal L^2(\Omega,\mu)$ has an orthonormal basis
consisting of all finite products
$$
e_j(x) =\langle j_0 , x_0\rangle \langle j_1, x_1 \rangle
\cdots \langle j_k, x_k \rangle \tag 8.13
$$
for $j=(j_0, j_1,\ldots,j_k,0,0,\ldots) \in \hat\Omega$, and
$x=(x_0,x_1,\ldots)\in \Omega$, where
$$
\langle j,x \rangle = \exp(2\pi i\, jx/n) \tag 8.14
$$
for $j,x\in \Bbb Z_n$.

We will now make a realization $T_1,\ldots, T_n$ of $\on$ on
$\Cal L^2(\Omega,\mu)$ which defines a shift without pure
normal invariant states. Any such realization has the form
$$
T_i = \sum_{j=1}^n S_j m_{ji}
$$
by \rom{(2.4)}, where $[m_{ji}]$ is a unitary matrix in $M_n
(\Cal B (\Cal L^2 (\Omega,\mu)))$. We take $[m_{ji}]$ to be
a diagonal matrix with $m_{ii}$ being the multiplication
operator on $\Cal L^2 (\Omega,\mu)$ defined by the function
$$
m_{ii}(x_0,x_1,x_2,\ldots)=\langle i,x_0 \rangle. \tag 8.15
$$

In formulas \rom{(2.5)} and \rom{(2.6)} above, we
introduced, for general pairs $(S_i),(T_i)$ in
$\rep(\on,\h)$, the {\it unitary transfer operator\/} $U$
which relates them. Recall that, for a general such pair,
$U$ is given by,
$$
U=\sum_j T_j S_j^* =\sum_i \sum_j S_i m_{ij} S_j^* ;
$$
and, for the present concrete pair, a calculation yields,
$$
(U\xi)(x_0,x_1,\ldots)=\langle  x_0,x_1 \rangle \xi
(x_0,x_1, \dots)\tag 8.16
$$
for $\forall \xi \in \Cal L^2 (\Omega,\mu)$ and for $\forall
x =(x_0,x_1, \ldots)\in \Omega$.

We are now ready to give the new shift associated with {\it
nearest neighbor states}. As we note in Remark 8.3 below,
this shift is {\it not\/} conjugate to any one of those from
Sections 6--7. Recall they were constructed from infinite
product states.

\proclaim{Theorem 8.1} Let  $(S_i)\in \rep(\on, \Cal
L^2(\mu))$ be given by \rom{(8.12)}, and let $\alpha$ be the
corresponding Haar shift. Let $T_i\in\rep(\on, \Cal
L^2(\mu))$ be given by, $T_i = S_i m_{ii}$, with the
functions $m_{ii}(\cdot)$ on $\Omega$ defined in
\rom{(8.15)}; and let, $\beta(A):= \sum_i T_i A T_i^*$, (for
$\forall A\in \Cal B(\Cal L^2(\mu))$) be the corresponding
endomorphism.

Then $\beta$ is a shift of Powers index $n$, and $\beta$
does not allow invariant vector states. The corresponding
state $\omega$ of $\uhf_n$ is given by
$$\align
\omega(e_{i_1 j_1}^{(1)} \otimes e_{i_2j_2}^{(2)} \otimes
\cdots \otimes e_{i_kj_k}^{(k)})
  & = \langle \1, T_{i_1} T_{i_2} \cdots T_{i_k} T_{j_k}^*
\cdots T_{j_1}^* \1 \rangle\\
  & = \frac{1}{n^k} \delta_{i_k j_k} \langle i_1, i_2
\rangle \langle i_2,i_3 \rangle
  \cdots \langle i_{k-1},i_k\rangle\\
&\quad \cdot\, \overline{\langle j_1, j_2\rangle}\,
\overline {\langle j_2,j_3 \rangle}
  \cdots \overline{\langle j_{k-1},j_k\rangle}.
\endalign$$\endproclaim

\demo{Proof} We have
$$
T_i^* = \bar m_{ii} S_i^*
$$
so by \rom{(8.11)} and \rom{(8.13)},
$$
(T_i^* \xi)(x_0,x_1,\ldots) = \overline{\langle i,x_0
\rangle} n^{-1/2} \xi
(i,x_0,x_1,\ldots) \tag 8.17
$$
for all $\xi \in \Cal L^2(\Omega,\mu)$. Assume now (ad
absurdum) that $\xi$ is a joint eigenvector of the
$T_i^*$'s:
$$
T_i^* \xi =\lambda_i \xi \tag 8.18
$$
where $\lambda_i \in \Bbb C$ and $\sum_{i=0}^{n-1}
|\lambda_i|^2=1$. Combining with \rom{(8.17)} we have
$$
\lambda_i \xi (x_0,x_1, x_2, \ldots) =\overline{\langle i,
x_0 \rangle} n^{-1/2} \xi (i,x_0,x_1, \ldots)
\tag 8.19
$$
for $i=0,\ldots,n-1$; that is,
$$
\xi(y_0,y_1,y_2,\ldots)
=\lambda_{y_0}n^{1/2}\langle y_0, y_1 \rangle \xi
(y_1,y_2,y_3,\ldots).  \tag 8.20
$$
By recursion,
$$\multline
\xi(y_0,y_1,y_2,\ldots)\\
=n^{m/2} \lambda_{y_0} \lambda_{y_1} \cdots \lambda_{y_{m-
1}}
\langle y_0, y_1 \rangle \langle y_1, y_2 \rangle \cdots
\langle y_{m-1},y_m \rangle \xi (y_m,y_{m+1},\ldots).
\endmultline\tag 8.21
$$
By the axiom of choice, there exist non-zero functions $\xi$
satisfying \rom{(8.19)}: One divides all strings
$(y_0,y_1,\ldots)$ into equivalence classes characterized by
having the same tail up to translations, and then one
assigns an arbitrary value of $\xi$ to one element in each
equivalence class and uses the recursion \rom{(8.19)} to
compute the value of $\xi$ on the other elements in the
class. We will now, however, argue that \rom{(8.19)} has no
nonzero solution $\xi\in \Cal L^2(\Omega,\mu)$. We will show
this by demonstrating that if $\xi\in \Cal L^2(\Omega,\mu)$
and $\xi$ satisfies \rom{(8.19)}, then $\xi$ is orthogonal
to all the vectors in the orthonormal basis \rom{(8.13)} for
$\Cal L^2(\Omega,\mu)$. The proof uses the Fourier transform
on the compact abelian group $\Omega$, and the corresponding
basis:  If $e_j(x)$ is the element given by \rom{(8.13)},
choose $m>k+1$ in \rom{(8.19)} to obtain
$$\multline
\widetilde{\xi}(j_0,j_1,\ldots,j_k,0,0,\ldots)  =: \langle
e_j,\xi \rangle
 =\int_\Omega \overline{e_j(y)} \xi(y)\,d \mu (y)\\
\shoveleft{=n^{-m/2} \sum_{y_0=1}^n\cdots \sum_{y_{m-1}=1}^n
   \lambda_{y_0} \lambda_{y_1}\cdots \lambda_{y_{m-1}}}\\
\shoveleft{\qquad \cdot  \overline{\langle j_0,y_0
\rangle}\, \overline{\langle j_1,y_1\rangle} \cdots
  \overline{\langle j_k,y_k \rangle} \cdot
 \langle y_0,y_1 \rangle \langle y_1,y_2 \rangle \cdots
  \langle y_{m-2}, y_{m-1} \rangle}\\
\shoveleft{\qquad\cdot \int_\Omega \langle y_{m-1},y_m
\rangle \xi (y_m, y_{m+1}, \ldots) \, d\mu
  (y_m,y_{m+1},\ldots)}\\
\shoveleft{= n^{-m/2} \sum_{y_0 =1}^n \cdots \sum_{y_{m -
1}=1}^n \lambda_{y_0} \lambda_{y_1} \cdots \lambda_{y_{m-1}}
\cdot
 \overline{\langle j_0, y_0 \rangle}\, \overline {\langle
j_1,y_1 \rangle}\,
  \cdots \overline {\langle j_k,y_k\rangle}}\\
\shoveleft{\qquad\cdot \langle y_0,y_1\rangle \langle
y_1,y_2 \rangle
   \cdots \langle y_{m-2},y_{m-1} \rangle \cdot
 \widetilde{\xi} (-y_{m-1},0,0,0,\ldots).\qquad\qquad\qquad}
\endmultline\tag 8.22$$

In the case $k=0$, $m=1$ an analogous formula takes the form
$$\align
\widetilde{\xi}(j_0,0,0,\ldots)
&=n^{-1/2} \sum_{y_0=1}^n \lambda_{y_0} \overline{\langle
j_0,y_0\rangle}
  \int \langle y_0,y_1 \rangle \xi (y_1,y_2, \ldots)\,
d\mu(y_1,y_2,\ldots)\tag 8.23\\
&=n^{-1/2} \sum_{y_0 =1}^n \lambda_{y_0} \overline{\langle
j_0,y_0 \rangle}\,
 \widetilde{\xi}(-y_0,0,\ldots).
\endalign$$
It follows, with $\widetilde{\xi}(j) =\widetilde{\xi}
(j,0,0,\ldots)$, that
$$\align
\sum_{j\in \Bbb Z_n} |\widetilde{\xi}(j)|^2
&= n^{-1} \sum_{j,y,z\in \Bbb Z_n} \overline{\lambda}_y
\lambda_z \langle j,y \rangle
  \overline{\langle j,z \rangle} \,
\overline{\widetilde{\xi} (-y)}\, \widetilde{\xi}(-z)\tag
8.24\\
& = n^{-1} \sum_{y,z\in \Bbb Z_n} n\delta (y-z)
\overline{\lambda}_y \lambda_z
\overline{\widetilde{\xi} (-y)}\, \widetilde{\xi} (-z)\\
&= \sum_{y\in \Bbb Z_n} |\lambda_y |^2\,|\widetilde{\xi} (-
y)|^2,
\endalign$$
so
$$
\sum_{j\in \Bbb Z_n} |\widetilde{\xi}(j)|^2
=\sum_{j\in\Bbb Z_n} |\lambda_j|^2\,|\widetilde{\xi}(-
j)|^2.\tag 8.25
$$
Since $\sum_{j\in \Bbb Z_n} |\lambda_j|^2 =1$, it follows
immediately that if at least two of the $\lambda_j$ are
nonzero, then $\widetilde{\xi}(j,0,0,\ldots) =0$ for all
$j\in\Bbb Z_n$. But the recursion relation \rom{(8.22)} then
implies that
$\widetilde{\xi}(j_0,j_1,\ldots,j_k,0,0,\ldots)=0$ for each
$(j_0,j_1,j_2,\ldots,j_k,0,0,\ldots)\in\hat\Omega$. It
follows that
$$
\xi =0\qquad\text{in $\Cal L^2(\Omega,\mu)$}\tag 8.26
$$
and \rom{(8.18)} has no nonzero solution.

If all $\lambda_j$ except one are zero, e.g.,
$(\lambda_j)=(1,0,\ldots,0)$, then it follows directly from
\rom{(8.21)} that
$$
\xi(y_0,y_1,y_2,\ldots)=0
$$
unless $(y_0,y_1,y_2,\ldots)=(0,0,0,\ldots)$. But this
single point has Haar measure zero, so again $\xi=0$ in
$\Cal L^2(\Omega,\mu)$.

This completes the proof that \rom{(8.18)} cannot have a
nontrivial solution. This means, by Theorem 4.1, that the
endomorphsim $\beta(A):=\sum_i T_i A T_i^*$, $A\in\bh$, {\it
cannot\/} have an invariant vector state.

To complete the proof of Theorem 8.1 we have to show that
$\beta$ really is a shift (using Theorem 3.3), and to
compute the corresponding state on $\uhf_n$. But using
\rom{(8.17)} and the corresponding formula
$$
(T_i\xi) (x_0,x_1,\ldots)
= \langle x_0,x_1 \rangle n^{1/2} \delta_{ix_0}
\xi(x_1,x_2,\ldots)
$$
one computes
$$\multline
(T_{i_1}\cdots T_{i_k} T_{j_k}^* \cdots T_{j_1}^*
\xi)(x_0,x_1,\ldots)\\
\shoveleft{= \delta_{i_1,x_0}\delta_{i_2,x_1} \cdots
\delta_{i_k,x_{k-1}}
   \langle x_0,x_1 \rangle \langle x_1,x_2 \rangle}\cdots
\langle x_{k-1},x_k \rangle\,\\
\shoveleft{\qquad\cdot\,\overline{\langle j_1,j_2\rangle}\,
  \overline{\langle j_2,j_3 \rangle}\,\cdots \,
\overline{\langle j_{k-1},j_k\rangle}\,
  \overline{\langle j_k,x_k \rangle}\,\xi
(j_1,\ldots,j_k,x_k,x_{k+1},\ldots)}\\
\shoveleft{=\delta_{i_1,x_0} \delta_{i_2,x_1}\cdots
\delta_{i_k,x_{k-1}}
  \langle i_1, i_2\rangle \langle i_2,i_3 \rangle
\cdots \langle i_{k-1},i_k\rangle \langle i_k,x_k \rangle
  \overline{\langle j_1,j_2 \rangle}\, \overline{\langle
j_2,j_3 \rangle}}\\
\shoveleft{\qquad \cdots \overline{\langle j_{k-
1},j_k\rangle}\, \overline{\langle j_k,x_k \rangle}\,
  \xi(j_1,\ldots,j_k,x_k,x_{k+1},\ldots).}
\endmultline$$
It follows immediately from this formula that the
representation of $\uhf_n$ on $\Cal L^2 (\Omega,\mu)$
defined by the $T_i$'s is irreducible, and thus by Theorem
3.3 $\beta$ is a shift. Furthermore
$$\multline
\omega (e_{i_1j_1}^{(1)} \otimes e_{i_2 j_2}^{(2)} \otimes
\cdots \otimes e_{i_k j_k}^{(k)})\\
\shoveleft{ = \langle \1, T_{i_1} \cdots T_{i_k} T_{j_k}^*
\cdots T_{j_1}^* \1\rangle}\\
\shoveleft{ = \frac{1}{n^{k+1}} \sum_{x_0} \cdots \sum_{x_k}
\delta _{i_1,x_0} \cdots
   \delta_{i_k,x_{k-1}}\cdot
 \langle i_1,i_2 \rangle \langle i_2,i_3 \rangle \cdots
   \langle i_{k-1},i_k \rangle \langle i_k,x_k \rangle}\\
\shoveleft{\quad\cdot\, \overline{\langle j_1,j_2\rangle}\,
\overline{\langle j_2,j_3\rangle}\,\cdots
   \overline{\langle j_{k-1},j_k \rangle}\,\overline{\langle
j_k,x_k\rangle}}\\
\shoveleft{= \frac{1}{n^{k+1}} \langle i_1,i_2 \rangle
\cdots \overline{\langle j_{k-1}, j_k \rangle}
   \sum_{x_k} \langle x_k,i_k -j_k \rangle}\\
\shoveleft{ = \frac{1}{n^k} \delta_{i_k j_k} \langle i_1,i_2
\rangle \langle i_2,i_3 \rangle
  \cdots \langle i_{k-1},i_k \rangle\cdot
  \overline{\langle j_1,j_2 \rangle}\, \overline{\langle
j_2,j_3 \rangle}
  \cdots \overline{\langle j_{k-1},j_k
\rangle.}\qquad\qquad\qquad}
\endmultline$$
This ends the proof of Theorem 8.1\qed\enddemo

\remark{Remarks 8.2} As already remarked after \rom{(8.21)},
the equation \rom{(8.19)} always has a continuum of function
solutions which are not measurable, and thus are not in
$\Cal L^2 (\Omega,\mu)$ or define states on $\on$ in any
reasonable sense. Also note that \rom{(8.19)} has the formal
infinite product ``solution''
$$
\xi(y_0,y_1,y_2,\ldots)=\prod_{k=0}^\infty
n^{1/2}\lambda_{y_k} \langle y_k, y_{k+1} \rangle.
$$
One way of stating Theorem 8.1 is that these infinite
products do not converge to a non-zero vector in $\Cal
L^2(\Omega,\mu)$. We will in the following remark on special
cases of \rom{(8.19)} where ``solutions'' exist which are
not in $\Cal L^2 (\Omega,\mu)$, and also supply some related
operator theoretic observations. Since, for the general case
of \rom{(8.18)}, or \rom{(7.10)} above, $\Cal L^2$-solutions
correspond to pure normal invariant states, the non $\Cal
L^2$ ``solutions'' correspond to states on $\on$ which are
{\it not normal\/} in the given Haar respresentation, and
the ``solutions'' give us a clue to what these states are,
namely the Cuntz states defined by the appropriate
$\lambda$'s. This lies at the heart of why one uses $C^*$-
algebras rather than merely Hilbert spaces in various
contexts: States which cannot be realized by vectors in the
Hilbert space, can be realized as states on an appropriate
$C^*$-algebra. In the analysis of quantum systems with
infinitely many degrees of freedom, examples of this abound
(sometimes under the name of the van Hove phenomenon); see
\cite{Br-Rob, p. 224}, \cite{Hov}, and \cite{Seg2}.

(i) In the special case where the vector $(\lambda_i)$ in
\rom{(8.18)} is $(1,0,\ldots,0)$, we noted that a possible
``eigenvector'' $\xi$ must then necessarily be a constant
times something like the delta mass at $0=(0,0,\ldots)$ in
the group $\Omega$. If $\xi$ shall define the appropriate
state on $\on$, it should rather be the ``square root'' of
the Dirac delta mass. This solution is not in $\Cal L^2
(\Omega,\mu)$, of course, unless the constant is zero.
Specifically, the assertion about $\xi$ in this special case
is that $\xi(x_0,x_1,\ldots)=0$ unless $x_0 =x_1 =\cdots
=0$.

(ii) The most interesting special case of \rom{(8.18)} turns
out to be the {\it equi-distribu\-tion}, $\lambda_i=n^{-
1/2}$ (for $\forall i$). In that case, the recursive formula
\rom{(8.20)} [for a possible $\Cal L^2$-solution $\xi$ to
\rom{(8.18)}] then takes the following geometric form: Using
\rom{(8.2)}, we may define the {\it isometric\/} operator
$S$ on $\Cal L^2(\Omega,\mu)$ by $S\xi:=\xi\circ\sigma$, and
\rom{(8.18)--(8.20)} become the single condition,
$$
\xi =US\xi \tag 8.27
$$
where $U$ is the unitary transfer operator from \rom{(2.6)}
and \rom{(8.16)}. \endremark

For this, moreover, the details for the \rom{(8.22)}
calculation simplify as follows. The present argument is
again based on the $\Omega$-$\Lambda$ duality and the
corresponding Fourier transform. Let $\lambda_i = n^{-1/2}$;
we supply the recursion. For $j\approx (j,0,0,\ldots )\in
\Lambda :=\hat\Omega$, we get:
$$\align
\widetilde{\xi}(j) &=\int_\Omega \overline{\langle x_0, j
\rangle}\, \xi (x_0,\ldots)\,d\mu(x_0,\ldots)\tag 8.28\\
&=n^{-2} \sum_{y_0} \sum_{y_1} \overline{\langle y_0,j
\rangle}\, \langle y_0,y_1 \rangle
 \int_\Omega \langle y_1, x_2 \rangle
\xi(x_2,\ldots)\,d\mu(x_2\cdots)\\
&=n^{-1} \sum_{y_1} \delta(j-y_1)\int_\Omega \langle y_1,x_2
\rangle \xi (x_2,\ldots) \,
 d\mu(x_2\cdots)\\
&=n^{-1} \int_\Omega \overline{\langle -j,x_2 \rangle}\, \xi
(x_2,\ldots)\,d\mu (x_2\cdots)\\
&=n^{-1} \widetilde{\xi} (-j)
\endalign$$
valid for $\forall j\in \Bbb Z_n :=\Bbb Z/n\Bbb Z$, and with
all summations being over $\Bbb Z_n$, again with $\Bbb Z_n$
viewed as an additive group. Replacing $j$ by $-j$  yields,
$\widetilde{\xi}(j)=0$, for $\forall j\in \Bbb Z_n$; or,
more specifically, $\widetilde{\xi}(j, 0, 0, \ldots)=0$, for
$\forall j \in \Bbb Z_n$, and \rom{(8.26)} follows as
before.

By a calculation quite analogous to the one above, we get,
$\forall (i_0, \ldots, i_k, 0,\ldots)\in \Lambda =\hat\Omega
=\coprod_0^\infty \Bbb Z_n$ , that
$$
\widetilde{\xi} (i_0,\ldots,i_k,0,0,\ldots)
=n^{-1} \overline{\langle i_0, i_1 \rangle}\,
\widetilde{\xi} (i_2 - i_0, i_3, \ldots, i_k, 0,0,\ldots).
\tag 8.29
$$
But then, by induction, $\widetilde{\xi}$ must vanish
identically on $\Lambda = \hat \Omega = \coprod_0^\infty
\Bbb Z_n$.

(iii) A more operator theoretic way to see that $US\xi=\xi$
implies $\xi=0$  is this: If $\xi \in \Cal L^2 (\Omega,
\mu)$ is arbitrary, one computes as above,
$$\align
\langle e_j, (US)^2\xi \rangle
=& {((US)^2 \xi)}\widetilde{\phantom{0}}\, (j_0, j_1, \ldots
j_k, 0, 0, \ldots)\\
= &n^{-1}\, \overline{\langle j_0,j_1 \rangle}\,
\widetilde{\xi} (j_2 - j_0, j_3, \ldots, j_k, 0,0, \ldots).
\endalign$$

Because of the $n^{-1}$ factor, it follows by iteration that
$$
|\langle e_j, (US)^{2m} \xi\rangle |\leq n^{-m}\|\xi\|.
$$

We will now show that \rom{(8.27)} has no solution by
proving that the unitary part of the {\it Wold
decomposition\/} of $T=US$ is zero. Recall from \cite{Nik}
that if $T$ is {\it any isometry\/} on $\h=\Cal
L^2(\Omega,\mu)$, i.e., $T^*T=1$, then $\h$ has a
decomposition, $\h=\h_1\oplus \h_2$ into $T$-invariant
subspaces such that $\left. T\right|_{\h_1}$ is unitary, and
$V=\left. T\right|_{\h_2}$ is a shift. That $V$ is a {\it
shift\/} means that $\lim_{n\rightarrow\infty} V^n \xi =0$
for any $\xi\in\h_2$. Put $\Cal L=\h_2 \ominus
V\h_2=\h\ominus T\h$. (If $(\xi_i)$ is an orthonormal basis
for $\Cal L$, and one defines $\xi_{ij} =V^j \xi_i =T^j
\xi_i$, then $(\xi_{ij})$ is an orthonormal basis for
$\h_2$. Thus $\h_2= \bigoplus_{m=0}^\infty V^m \Cal L$, and
$V$ decomposes into a direct sum of $\dim \Cal L$
copies of  the {\it Hilbert shift}, defined by, $\xi_{ij}
\mapsto \xi_{i,j+1}$, for fixed $i$, and $j=0,1,\ldots$;
$\dim \Cal L$ is called the {\it multiplicity\/} of the
shift.)  The subspaces $\h_1$ and $\h_2$ may be identified
through the formuli
$$
\h_1=\bigcap_m T^m \h =\bigcap_m \{\xi \in\h :
\|T^{*m}\xi\|= \|\xi\|\}
$$
and
$$
\h_2 =\h_1^\bot =\bigoplus_{m=0}^\infty T^m \Cal L.
$$

Returning to our specific isometry $T=US$, we have to show
that $\h_1 =0$: Let $\xi\in \h_1$. Then $\xi\in T^{2m}\h$,
so for each $m$ there exists a $\xi_m \in \h$ with
$\xi=T^{2m}\xi_m$. But then
$$
\|\xi_m\| = \|T^{2m}\xi_m \| =\|\xi \|,
$$
and hence
$$\align
|\langle e_j,\xi \rangle | &=|\langle e_j,
T^{2m}\xi_m\rangle| \\
&\leq n^{-m} \|\xi_m\|=n^{-m}\|\xi \|.
\endalign$$
Letting $m\rightarrow \infty$, we see that
$$
\langle e_j, \xi \rangle =0,
$$
and, since $j\in \Lambda$ is arbitrary, the desired
conclusion, $\xi=0$, follows. We conclude that $T=US$ is a
{\it shift\/} on $\Cal L^2(\mu)$, equivalently a completely
non-unitary isometry. This seems of independent interest as
the isometry $S$ is {\it not\/} a shift, recall $S^*\1=\1$.
An inspection also reveals that the shift $T$ has the
multiplicity $(n-1)\cdot\infty$ where $n$ is the index of
the Haar shift. Of course, the infinite product,
$$\multline
\xi(x_0,x_1,x_2,\ldots) := \prod_{k=0}^\infty \langle x_k,
x_{k+1} \rangle\\
\shoveleft{\text{(or, more formally,}\prod_{k=0}^\infty e^{i
2\pi x_k x_{k+1} /n} = e^{i(2\pi/n)\sum_{k=0}^\infty x_k
x_{k+1}})\qquad\qquad\qquad\qquad}
\endmultline$$
is a ``formal'' solution to \rom{(8.27)}; but our present
considerations imply that this infinite product is indeed
purely formal, and {\it not\/} convergent in $\Cal
L^2(\Omega,\mu)$. Specifically (ad absurdum), convergence in
$\Cal L^2(\mu)$ would put the limit-function, $\xi(x)$ (for
$x\in \Omega$), in $\bigcap_{k=1}^\infty T^k (\Cal
L^2(\mu))$. But this intersection {\it is\/} the unitary
term in the Wold-decomposition, and we proved that it must
be zero.

Note furthermore that our stronger conclusion, based on this
Wold decomposition argument, is the assertion that there can
be no sequence $(\xi_k)_{k=1}^\infty$ in $\Cal L^2(\mu)$
such that the limit, $\lim_{k\rightarrow \infty} T^k \xi_k$,
exists in $\Cal L^2(\mu)$ and is non-zero.

(iv) With the notation
$$
e_{ij} := e_{i_1 j_1}^{(1)} \otimes \cdots \otimes
e_{i_kj_k}^{(k)} \otimes 1 \otimes \cdots \in \uhf_n,
$$
and
$$
i=(i_1,\ldots,i_k,0,0,\ldots)\in \Lambda = \coprod_1^\infty
\Bbb Z_n = \widehat{\left(\prod_1^\infty \Bbb Z_n \right)},
$$
the formula for the state $\omega$ in Theorem 8.1 above is
$$
\omega(e_{ij}) =n^{-k} \delta(i_k -
j_k)\xi(i)\overline{\xi(j)}
$$
where the function $\xi(\cdot)$ is defined (as in (iii)
above) on $\Lambda$ by
$$
\xi(i)=\prod_{p=1}^\infty \langle i_p,i_{p+1}\rangle,\tag
8.30
$$
and, for positive definite functions $F(\cdot,\cdot)$ on
$\Lambda\times\Lambda$, there are corresponding states
$\omega_F$ on $\uhf_n$ given by the more general formula
$$
\omega_F(e_{ij}) := F(i,j)\xi(i)\overline{\xi(j)}\tag 8.31
$$
for $(i,j)\in \Lambda\times\Lambda$ having the same length.
When $F$ is so chosen, the object is to identify operators
$T_i$, depending on $F$, and satisfying the Cuntz-relations,
such that
$$\omega_F(e_{ij}) =
\langle \1,T_{i_1}\cdots T_{i_k}T_{j_k}^* \cdots T_{j_1}^*
\1\rangle_{\Cal L^2}
$$
is given by the expression on the right hand side in
\rom{(8.31)} and $\1$ denotes the constant unit function on
$\Omega$. Specifically we may get such states $\omega_F$ in
the set $P$ from Section 5 as follows: Let $\omega =
\bigotimes_1^\infty \omega_k$ (each $\omega_k$ is a state on
$M_n$) be a product state in $P$ as described in Example
5.5, and for $i=(i_1,\ldots, i_k,0,0,\ldots)$,
$j=(j_1,\ldots,j_k,0,0,\ldots)$ in $\Lambda$, let
$$
F_\omega(i,j) := \prod_{m=1}^k \omega_m(e_{i_m j_m}^{(m)})
\cdot \sum_{r\in \Bbb Z_n} \langle r,i_k -j_k \rangle
\omega_{k+1} (e_{rr}^{(k+1)}).
$$
Then it can be shown that the corresponding state
$\omega_{F_\omega}$ in \rom{(8.31)} is in $P$ (details in a
sequel paper), and we get an associated element in
$\rep_s(\on,\Cal L^2)$. Furthermore, we may choose the
product state $\bigotimes_k \omega_k$ in $P$ such that the
corresponding {\it shift\/} $\beta$ on $\Cal B(\Cal L^2)$,
i.e., $\beta(A)=\sum_i T_i AT_i^*$, is {\it non-conjugate\/}
to the one from Theorem 8.1, and also not to those from
Sections 6--7. Note that the function $\xi(\cdot)$ in
\rom{(8.31)} is well defined on the subgroup $\Lambda$ of
$\Omega$, but, as noted in (iii) above, it is not
sufficiently {\it almost periodic\/} to extend naturally to
the compactification $\Omega$.

\remark{Remark 8.3} It can be proved that if $\omega$ is the
state on $\uhf_n$ defined in Theorem 8.1, and $\omega'$ is
any infinite tensor product of pure states on $\uhf_n$, then
$$
\|\omega\circ\sigma^m - \omega' \circ \sigma^m \| =2
$$
for all $m\in \Bbb N$. If $\omega'\in P$, it follows from
Lemma 5.4 that the corresponding shifts of $\bh$ are non-
conjugate, and hence the shift considered in this section is
not conjugate to any one of those discussed in Sections 6
and 7. The proof will be deferred to a subsequent paper
where nearest neighbor states and similar states will be
treated more systematically.
\endremark

\head 9. Extending Unital Endomorphisms to Automorphisms
\endhead

In \cite{Arv2}, \cite{AK} it was proved that a continuous
one-parameter semigroup of unital $^*$-endomorphisms of
$\bh$ has an extension to a group of $^*$-automorphisms of
$\Cal B(\h \otimes\h)$ when $\bh$ is embedded as $1
\otimes\bh$. Using techniques from \cite{PR}, let us
establish a similar (but simpler) result for a single
endomorphism.

\proclaim{Theorem 9.1} Let $\alpha$ be a unital $^*$-
endomorphism of $\bh$ of index $n$, and embed $\bh$ into
$\Cal B(\h \otimes \h)$ as $1 \otimes \bh$. Then $\alpha$
has an extension $\beta$ to $\Cal B(\h \otimes \h)$ such
that $\beta$ is a $^*$-automorphism. Furthermore, if
$\alpha$ is a shift, and $H_0$ is the Hilbert space of
dimension $n$, and $M_n=\Cal B(H_0)$, then $\Cal B(\h
\otimes \h)$ contains $M_{n^\infty} = \bigotimes _{k=-
\infty}^\infty M_n$ as a weakly dense $C^*$-subalgebra in
such a manner that the restriction of $\beta$ to
$M_{n^\infty}$ is just the two-sided right shift, and
$\bigotimes_{k=0}^\infty M_n \subseteq \bigotimes_{k=-
\infty}^\infty M_n$ is weakly dense in $1\otimes
\bh$.\endproclaim

\demo{Proof} Since $\alpha(\bh)' \cong \Cal B(H_n)$, we have
a tensor product decomposition
$$
\h = H_n \otimes \Cal K
$$
such that
$$
\alpha(\bh) = 1_{H_n} \otimes \Cal B(\Cal K).
$$
Let $\alpha' : \bh \rightarrow \Cal B(\Cal K)$ be the
corresponding $^*$-isomorphism such that
$$
\alpha(A) = 1 \otimes \alpha' (A)
$$
for $A\in \bh$. Choose a  particular unit vector in $H_0$,
and let $H' = \bigotimes_{k=-\infty}^{-1} H_0$ be the
corresponding von Neumann reduced tensor product. (For the
first part of Theorem 9.1, we do not need any structure on
$H'$ other than it is a separable infinite dimensional
Hilbert space.) First, let $\beta'$ be {\it any\/} $^*$-
isomorphism $\Cal B(H') \rightarrow \Cal B(H' \otimes H_0)$
and define
$$
\beta: \Cal B(H' \otimes \h) \rightarrow \Cal B(H' \otimes
\h)
$$
by
$$
\beta (B\otimes A) = \beta' (B) \otimes \alpha'(A)
$$
for $B\in \Cal B(H')$, $A\in \bh$ and the {\it last\/}
tensor product is according to the decomposition
$$
H' \otimes \h = (H' \otimes H_0)\otimes \Cal K.
$$
Then $\beta$ is a $^*$-automorphism extending $\alpha$. For
the last part of the theorem we define $\beta'$, more
specifically, as the $^*$-isomorphism, $\Cal B(\bigotimes_{-
\infty}^{-1} H_0) \rightarrow \Cal B(\bigotimes_{-\infty}^0
H_0)$, implemented by the right-shift, $U:\bigotimes_{-
\infty}^{-1} H_0 \rightarrow \bigotimes_{-\infty}^0 H_0$,
defined by
$$
U(\cdots \psi_{-3} \otimes \psi_{-2} \otimes \psi_{-1} )
= \cdots \underset {-2} \to {\psi_{-3}} \otimes \underset{-
1} \to {\psi_{-2}} \otimes
  \underset 0 \to {\psi_{-1}}.
$$
Now, if $N_0 = \alpha(\bh)' \cap \bh$ and, inductively
$$
N_{k+1} = \alpha (N_k),\qquad k=0,1,\ldots
$$
then by \cite{Pow2, Lemma 2.1}, or the ideas surrounding
\rom{(3.6)} in the proof of Theorem 3.3, it follows that the
$N_k$'s are mutually commuting $I_n$ factors, with
$\{\bigcup_k N_k\}' = \Bbb C(1)$. Putting $N_{-k}$ equal to
the bounded operators of the $-k'$th tensor factor in
$\bigotimes_{-\infty}^{-1} H_0$ (tensor $1$ on the remaining
factors), it follows that all the $N_k$'s mutually commute,
$\beta(N_k) = N_{k+1}$ for $k\in \Bbb Z$, and the $C^*$-
algebra generated by the $N_k$'s is weakly dense in $\Cal
B(H' \otimes \h)$. Theorem 9.1 follows. \qed\enddemo

\remark{Remark 9.2} If $\alpha$ is a shift, and $\beta$ and
$\Cal B(H'\otimes \h)$ are constructed according to the
recipe above, then all elements in the weakly dense $^*$-
subalgebra $\bigcup_{k=-1}^{-
\infty}\left(\bigotimes_k^\infty M_n\right)$ of $\Cal B(H'
\otimes \h)$ will ultimately be mapped into $1_{H'} \otimes
\bh$. Thus any asymptotic property of $\alpha$ (such as
having an absorbing state), readily translates into a
similar property for $\beta$.
\endremark

\subhead Acknowledgments \endsubhead The present paper came
about in connection with the second named author's (Palle
Jorgensen) visit to the University of Oslo, with support
from {\it the Norwegian Research Council}, and the visit by
the last two authors (Palle Jorgensen and Geoffrey L. Price)
to the Massachusetts Institute of Technology in connection
with {\it the 1994 von Neumann Symposium}, and with support
from the {\it U.S. National Science Foundation}. Hospitality
from the respective hosts (O. Bratteli and E. St\o rmer in
Oslo, and I.E. Segal at MIT) is greatly appreciated. In
addition to the two funding agencies already mentioned, this
work was also supported by a {\it U.S.--Western Europe NSF
grant\/}, by  a {\it NATO grant}, and (for Palle Jorgensen)
a {\it University of Iowa Faculty Scholarship Award}. We
acknowledge helpful conversations with G.A.~Elliott,
R.T.~Powers, E.~St\o rmer, and S.~Pedersen. The paper was
finished while Ola Bratteli visited {\it The Fields
Institute for Research in Mathematical Sciences\/} with
Norwegian Research Council funding. He is grateful for
patience from the host G.A.~Elliott there. Finally we wish
to thank M.~Laca for generously making many valuable
suggestions to a preliminary version of our paper, and for
allowing us to quote from \cite{Lac}.

\enddocument